\documentclass[hidelinks]{siamart190516}

\usepackage{smalljeff}
\usepackage{amsmath,amssymb,mathtools}
\usepackage{booktabs,enumitem}
\usepgfplotslibrary{fillbetween, groupplots, patchplots}
\usetikzlibrary{patterns}
\usepackage{datatool}
\usepackage{xfrac}

\usepgfplotslibrary{external} 
\newcommand{\linelabel}[1]{}

\usepackage{makecell}

\newcommand{\TheTitle}{A Lipschitz Matrix for Parameter Reduction in Computational Science}
\newcommand{\TheAuthors}{Jeffrey M. Hokanson and Paul G. Constantine}
\headers{Lipschitz Matrix}{\TheAuthors}

\DeclareMathOperator{\vol}{{\mathrm vol}}

\newcommand{\fmin}{f_{\min}}
\newcommand{\fmax}{f_{\max}}
\newcommand{\fmid}{f_{\mathrm{mid}}}

\title{{\TheTitle}\thanks{Submitted to the editors 3 June 2019.
\funding{This work is supported by DARPA's program Enabling Quantification of Uncertainty in Physical Systems
(EQUiPS). 
}}}
\author{Jeffrey M. Hokanson\thanks{
	Department of Computer Science, 
	University of Colorado Boulder,
	1111 Engineering Dr, Boulder, CO 80309,
	(\email{Jeffrey.Hokanson@colorado.edu}, \email{Paul.Constantine@colorado.edu}).
	}
	\and Paul G. Constantine\footnotemark[2]}

\begin{document}
\maketitle
\begin{abstract}
We introduce the Lipschitz matrix:
a generalization of the scalar Lipschitz constant
for functions with many inputs.
Among the Lipschitz matrices compatible a particular function,
we choose the smallest such matrix in the Frobenius norm
to encode the structure of this function.
The Lipschitz matrix then provides a function-dependent metric on the input space.
Altering this metric to reflect a particular function
improves the performance of many tasks in computational science.
Compared to the Lipschitz constant, 
the Lipschitz matrix reduces the worst-case cost of approximation, integration, and optimization;
if the Lipschitz matrix is low-rank,
this cost no longer depends on the dimension of the input,
but instead on the rank of the Lipschitz matrix
defeating the curse of dimensionality.
Both the Lipschitz constant and matrix 
define uncertainty away from point queries of the function
and by using the Lipschitz matrix we can reduce uncertainty.
If we build a minimax space-filling design of experiments in the Lipschitz matrix metric,
we can further reduce this uncertainty.
When the Lipschitz matrix is approximately low-rank,
we can perform parameter reduction by constructing
a ridge approximation whose
active subspace is the span of the dominant eigenvectors of the Lipschitz matrix.
In summary, the Lipschitz matrix provides a new tool
for analyzing and performing parameter reduction
in complex models arising in computational science.
\end{abstract}
\begin{DOI}
\end{DOI}
 \section{Introduction}
With the increasing sophistication of computer models, 
practitioners in science and engineering often confront the \emph{curse of dimensionality}---the phenomena that for many relevant computational tasks, 
obtaining the desired solution requires a computational cost that grows exponentially 
with the number of input parameters~\cite{Don00,TW98}.
To mitigate this curse,
practitioners may introduce a lower dimensional reparameterization 
of the model inputs  
yielding a similar output to the original.
When applicable, this parameter reduction enables otherwise high-dimensional 
computations to exploit low-dimensional constructed parameters.
Standard engineering practice uses a \emph{global sensitivity analysis}~\cite{SRA+08}
to identify a subset of input parameters that have relatively little affect on the output.
Fixing these insignificant parameters at a nominal value
provides a lower-dimensional parameterization using the remaining parameters~\cite{STG+07}.
Another approach reparameterizes the model
using a few linear combinations of the input variables
defining an \emph{active subspace} along which the model varies~\cite{Con15}. Active subspaces include subset-based approaches
since any subset can be encoded as the span of columns of the identity matrix.
More generally, a nonlinear reparameterization of the input-output map can be used~\cite{HKB+19}.
Here we introduce the \emph{Lipschitz matrix}---a generalization of the scalar Lipschitz constant.
The Lipschitz matrix can identify an active subspace,  
motivates a design of experiments, 
defines uncertainty away from function evaluations,
and yields improved error bounds that can mitigate the curse of dimensionality.

\subsection{Definition\label{sec:intro:def}}
We define the Lipschitz matrix analogously to the Lipschitz constant.
Given a domain $\set D \subset \R^m$,
the scalar Lipschitz constant $L \in \R_+$ defines a set of  
\emph{scalar Lipschitz functions} denoted $\set L(\set D, L)$ that satisfy
\begin{equation}\label{eq:scalar_lipschitz}
	\set L(\set D, L) \coloneqq
		\lbrace f:\set D\to \R \ : \  |f(\ve x_1) - f(\ve x_2)| \le L \|\ve x_1 - \ve x_2\|_2,
		\ \ve x_1, \ve x_2 \in \set D\rbrace.
\end{equation}
The \emph{Lipschitz matrix} changes this definition,
moving the constant $L$ inside the norm and promoting it to a matrix $\ma L \in \R^{m\times m}$.
This defines \emph{matrix Lipschitz functions} 
\begin{equation}\label{eq:matrix_lipschitz}
	\set L(\set D, \ma L) \coloneqq
		\lbrace f:\set D\to \R \ : \  |f(\ve x_1) - f(\ve x_2)| \le \|\ma L(\ve x_1 - \ve x_2)\|_2,
		\  \ve x_1, \ve x_2 \in \set D\rbrace.
\end{equation}
In our notation, the type of the second argument of $\set L(\cdot, \cdot)$
indicates whether this set refers to the scalar \cref{eq:scalar_lipschitz}
or matrix \cref{eq:matrix_lipschitz} case.
Additionally, we define \emph{$\epsilon$-Lipschitz functions}
that are near by Lipschitz functions:
\begin{equation}\label{eq:Lmat_epsilon}
	\set L_\epsilon(\set D, \ma L) \coloneqq 
		\lbrace f: \set D\to \R: 
			|f(\ve x_1) - f(\ve x_2)| \le \|\ma L(\ve x_1 - \ve x_2)\|_2 + \epsilon,
			\  \ve x_1, \ve x_2 \in \set D 
		\rbrace.
\end{equation}
This function class is useful when 
analyzing functions with \emph{computational noise}~\cite{MW11}
or irrelevant oscillations;
see examples in \cref{sec:epsilon}.

\subsection{Equivalence}
We refer to both scalar~\cref{eq:scalar_lipschitz} and matrix~\cref{eq:matrix_lipschitz} 
Lipschitz functions as simply \emph{Lipschitz functions} as these two sets are nested.
All scalar Lipschitz functions $f \in \set L(\set D, L)$
are also matrix Lipschitz functions with $f\in \set L(\set D, L \ma I)$:
\begin{equation}
	|f(\ve x_1) - f(\ve x_2)| \le L \| \ve x_1 - \ve x_2\|_2 = \| (L\ma I)(\ve x_1 - \ve x_2)\|_2.
\end{equation}
Hence $\set L(\set D, L) = \set L(\set D, L \ma I)$.
Similarly, all matrix Lipschitz functions $f\in \set L(\set D, \ma L)$ are also scalar Lipschitz functions 
with $f\in \set L(\set D, \|\ma L\|_2)$ as 
\begin{equation}
	|f(\ve x_1) - f(\ve x_2)| \le \|\ma L(\ve x_1 - \ve x_2)\|_2 \le \|\ma L\|_2 \|\ve x_1 - \ve x_2\|_2.
\end{equation}
Hence $\set L(\set D, \ma L) \subseteq \set L(\set D, \|\ma L\|_2)$.
This nesting also holds for $\epsilon$-Lipschitz functions.

\subsection{Smallest Lipschitz Matrix\label{sec:intro:smallest}}
Given a particular function, 
we seek the smallest possible Lipschitz matrix
to tighten our results.
The challenge using the Lipschitz matrix
is there is no natural ordering unlike the scalar Lipschitz constant.
Hence we must impose an ordering.

For the scalar Lipschitz constant,
the ordering of the real line provides a natural ordering
for Lipschitz constants.
The smallest Lipschitz constant for a function $f$ is
\begin{equation}
	\min_{L\in \R_+} \ L \quad \text{subject to} \quad
	|f(\ve x_1) - f(\ve x_2)| \le L \|\ve x_1 - \ve x_2\|_2 \quad  \ve x_1, \ve x_2 \in \set D.
\end{equation}
Unlike scalars, matrices have no natural ordering.
However by invoking the polar decomposition~\cite[Thm.~7.3.1]{HJ13}
we can define a partial order for Lipschitz matrices.
Any matrix $\ma L\in \R^{m\times m}$ has a polar decomposition into 
the product $\ma L = \ma Q \ma P$ where $\ma Q \in \R^{m\times m}$
has orthonormal columns and $\ma P \in \R^{m\times m}$ is a symmetric positive semidefinite matrix,
denoted $\ma P \in \mathbb{S}_+^{m \times m}$.
As the 2-norm is unitarily invariant,
\begin{equation}\label{eq:lipschitz_invariant}
	\|\ma L(\ve x_1 - \ve x_2)\|_2 = \|\ma Q \ma P(\ve x_1 - \ve x_2)\|_2 = \|\ma P(\ve x_1 - \ve x_2)\|_2.
\end{equation}
Thus without loss of generality, we can assume a Lipschitz matrix $\ma L$ 
is symmetric positive semidefinite.
Positive semidefinite matrices have a natural partial ordering:
the \emph{Loewner partial order}~\cite[Def.~7.7.1]{HJ13}
where for $\ma A, \ma B \in \mathbb S_+^{m\times m}$,
we write $\ma A \preceq \ma B$ if $\ma B - \ma A$ is positive semidefinite.
This is only a partial order as when $\ma B - \ma A$ is indefinite, 
we cannot order $\ma A$ and $\ma B$.

To define the smallest Lipschitz matrix, 
we choose a total order compatible with the partial order.
Two convenient choices are the trace and determinant of $\ma L$~\cite[Cor.~7.7.4d,e]{HJ13}.\linelabel{line:HoJo}
In the majority of our results, 
we use the squared Frobenius norm of $\ma L$,
i.e., the squared sum of the eigenvalues of $\ma L$ which obeys the partial order by~\cite[Cor.~7.7.4c]{HJ13},
as this yields an amenable optimization problem as discussed in~\cref{sec:estimate}.

\subsection{Derivatives\label{sec:intro:der}}
The Lipschitz matrix bounds the derivatives of $f$.
Suppose $f$ is differentiable at $\ve x$ in the interior of $\set D$.
Then for any $\ve p \in \R^m$ there is some $\delta > 0$
such that $\ve x +\delta \ve p \in \set D$.
From the Lipschitz matrix constraint~\cref{eq:matrix_lipschitz},
\begin{equation}
	|f(\ve x+\delta \ve p) - f(\ve x)|^2 \le \|\ma L(\ve x +\delta \ve p - \ve x)\|_2^2
		= \delta^2 \ve p^\trans \ma L^\trans \ma L \ve p.
\end{equation}
Dividing by $\delta^2$ and taking the limit as $\delta \to 0$,
we bound the gradient of $f$, $\nabla f$, by
\begin{equation}
	(\ve p^\trans \nabla f(\ve x))^2 = \ve p^\trans \nabla f(\ve x) \nabla f(\ve x)^\trans \ve p
		\le \ve p^\trans \ma L^\trans \ma L \ve p
\end{equation}
As this holds for all $\ve p$,
we write this compactly using the Loewner partial order:
\begin{equation}\label{eq:grad_order}
	\nabla f(\ve x) \nabla f(\ve x)^\trans \preceq \ma L^\trans \ma L.
\end{equation}
If $f$ is not differentiable,
a similar result holds for any Gateaux derivative.
Note, due to the presence of $\epsilon>0$ in the definition of $\epsilon$-Lipschitz functions~\cref{eq:Lmat_epsilon},
the derivatives of $f$ do not provide a lower bound on the $\epsilon$-Lipschitz matrix.\linelabel{line:noder}

\subsection{Connection to ridge functions}
A \emph{ridge function}~\cite{Pin15} depends only on a few linear combinations of its input variables;
i.e., a function $f:\R^m\to \R$ is a ridge function if there is a function $g$ such that
\begin{equation}
	f(\ve x) = g(\ma U^\trans \ve x) 
	\quad \text{where} \quad 
	g: \R^n\to \R, \quad 
	\ma U \in \R^{m\times n}, \quad
	 m > n.
\end{equation}
We call $g$ the \emph{ridge profile}, $n$ the \emph{ridge dimension},
and the range of $\ma U$ the \emph{active subspace}.
The Lipschitz matrix is intimately connected to ridge functions:
informally, $f$ is a ridge function if and only if it has a low-rank Lipschitz matrix.
The following theorem makes this precise.
\begin{theorem}\label{thm:ridge}
	Suppose $f: \set D \to \R$ is a Lipschitz function with a Lipschitz constant $L$.
	Then $f$ is a ridge function with ridge dimension $n$ 
	if and only if there exists a $\ma L \in \R^{m\times m}$ 
	with rank $n$ such that $ f \in \set L(\set D, \ma L)$.
\end{theorem}
\begin{proof}
	Suppose $f \in \set L(\set D, \ma L)$ where $\ma L$ is rank $n$.
	Let $\ma U \in \R^{m\times n}$ be an orthonormal basis for the range of $\ma L$.
	If $\ve x_1, \ve x_2 \in \set D$ with $\ma U^\trans(\ve x_1 - \ve x_2) = 0$, then
	\begin{equation}
		|f(\ve x_1) - f(\ve x_2)| \le \|\ma L(\ve x_1 - \ve x_2)\|_2 = 0
	\end{equation}	
	as $\ve x_1 - \ve x_2$ is in the nullspace of $\ma L$.
	Hence $f$ is constant in all directions in the nullspace of $\ma U$
	and thus there is a $g: \R^n \to \R$ such that $f(\ve x) = g(\ma U^\trans \ve x)$.

	Suppose $f(\ve x) = g(\ma U^\trans \ve x)$
	where $\ma U\in \R^{m\times n}$ has orthonormal columns.
	As $f$ is Lipschitz, so too is $g$.
	Let $\ma L_g\in \mathbb{S}_+^{n\times n}$ be a full-rank Lipschitz matrix for $g$,
	then
	\begin{equation}
		\begin{split}
		|f(\ve x_1) - f(\ve x_2)| &
		= |g(\ma U^\trans \ve x_1) - g(\ma U^\trans \ve x_2)| \\
		&\le \| \ma L_g[ (\ma U^\trans \ve x_1) - (\ma U^\trans \ve x_2)]\|_2 
		= \| \ma U\ma L_g \ma U^\trans(\ve x_1 - \ve x_2)\|_2.
		\end{split}
	\end{equation}
	Hence $f\in \set L(\set D, \ma U\ma L_g \ma U^\trans)$
	and $\ma U\ma L_g\ma U^\trans$ is rank $n$.
\end{proof}

Thus if we identify a low-rank Lipschitz matrix for a function,
the range of the Lipschitz matrix defines an active subspace.
However in our experience exact ridge functions are rare.
More frequently, a function will be approximately a ridge function;
we informally call these functions \emph{ridge-like}.
In the context of Lipschitz matrices,
we say $f$ is ridge-like if it has an approximately low-rank Lipschitz matrix;
then the dominant eigenspace defines an active subspace (see \cref{thm:error_bound}).
Other approaches for identifying and approximating ridge-like functions include:
polynomial ridge approximation~\cite{CEHW17, HC18},
\emph{sufficient dimension reduction} techniques from statistical regression~\cite{glaws2019},
and the mean gradient outer-product (MeGO)~\cite{Con15}.

\subsection{Comparison to MeGO}
Given a probability measure $\mu$,
the mean gradient outer-product (MeGO)~\cite[eq.~(3.2)]{Con15} is
\begin{equation}\label{eq:EGO}
	\ma C \coloneqq \mathbb{E} [ \nabla f \nabla f^\trans] = 
		\int_{\ve x\in \set D} \nabla f(\ve x) \nabla f(\ve x)^\trans \D \mu(\ve x)
		\in \R^{m\times m}.
\end{equation} 
The dominant eigenspaces of $\ma C$ provide one way 
to identify the active subspace of $f$.

When $f$ is differentiable almost everywhere on its domain,
we can bound the MeGO matrix by the Lipschitz matrix.
From~\cref{eq:grad_order}, the Lipschitz matrix bounds the gradient:
$\nabla f(\ve x) \nabla f(\ve x)^\trans \preceq \ma L^\trans \ma L$.
Then as $\mu$ is a probability measure
\begin{equation}\label{eq:C_L_bound}
	\ma C = \int_{\ve x\in \set D} \nabla f(\ve x) \nabla f(\ve x)^\trans \D \mu(\ve x)
	\preceq \ma L^\trans \ma L.
\end{equation}
This bound is tight when $f$ is a linear function $f(\ve x) = \ve a^\trans \ve x $
in which case $\ma C  = \ve a \ve a^\trans$,
$\ma L = \|\ve a\|_2^{-1} \ve a \ve a^\trans$,
and $\ma L^\trans\ma L = \ve a \ve a^\trans$.
Like the Lipschitz matrix, 
dominant eigenspace of $\ma C$ 
identifies an active subspace~\cite[Thm.~1]{Con15}.
As a corollary,
if $f$ is a ridge function,
then the nullspaces of $\ma C$ and $\ma L$ are the same.

\subsection{Applications of the Lipschitz Matrix}
The Lipschitz matrix provides
both analytical tools and practical results.
As we show in \cref{sec:ibc}, 
replacing the Lipschitz constant with the Lipschitz matrix
allows us to approximate, integrate, and optimize
functions using fewer evaluations.
If we can identify a low-rank Lipschitz matrix,
then the cost of these tasks no longer scales with the number of parameters,
but instead the rank of the Lipschitz matrix.
In practice, there are few functions
for which we can compute the Lipschitz matrix exactly.
We show in \cref{sec:estimate} that we can estimate the Lipschitz matrix
from arbitrary combinations of function evaluations $f(\ve x_i)$ and gradients $\nabla f(\ve x_k')$
by solving a semidefinite program.
We can use the estimated Lipschitz matrix to then
identify an active subspace (\cref{sec:ridge}),
guide the design of experiments (\cref{sec:design}),
construct ridge approximations for dimension reduction (\cref{sec:dimension}),
and quantify uncertainty (\cref{sec:uq}).

\subsection{Reproducibility}
Following the principles of reproducible research, 
we provide code implementing the algorithms described in this paper and 
scripts generating the data appearing in the figures and tables
available at {\tt \url{http://github.com/jeffrey-hokanson/PSDR/}}.

 \section{Algorithm Complexity\label{sec:ibc}}
In this section we use results from \emph{information-based complexity}~\cite{TW98, TWW88}
to bound the worst-case optimal cost of approximation, integration, and optimization
for Lipschitz matrix functions.
These results parallel similar results for scalar Lipschitz functions~\cite{Suk92}
with one important distinction.
For scalar Lipschitz functions on $\set D \subset \R^m$,
each of these three tasks requires $\order(\epsilon^{-m})$ function evaluations
to obtain $\epsilon>0$ accuracy---this exponential growth in dimension is the \emph{curse of dimensionality}~\cite{Don00}.
For functions with a rank-$r$ Lipschitz matrix
these tasks require only $\order(\epsilon^{-r})$ function evaluations.
If $r< m$, then complexity is independent of dimension
and we have mitigated the curse of dimensionality.
The key ingredient is using the Lipschitz matrix
to provide a (pseudo-)metric on the domain 
\begin{equation}
	d_{\ma L}(\ve x, \ve y) \coloneqq \|\ma L(\ve x - \ve y)\|_2
	\qquad 
	\ve x, \ve y\in \set D.
\end{equation}
Using this metric, we show the complexity of approximation, integration, and optimization
is proportional to \emph{$\epsilon$-internal covering number} of $\set D$:
\begin{equation}\label{eq:cover}
	N_\epsilon(\set D;\ma L)  \coloneqq  \argmin_{M, \lbrace \ve x_j\rbrace_{j=1}^M \subset \set D}   M 
	\ \  \text{such that}\  \ 
	\set D \subseteq \bigcup_{j=1}^M \set B_\epsilon(\ve x_j;\ma L),
\end{equation}
where $\set B_{\epsilon} (\ve x_j;\ma L) \coloneqq \lbrace \ve x : \|\ma L(\ve x - \ve x_j)\|_2 \le \epsilon\rbrace$\linelabel{line:ball}
is the $\epsilon$-ball in $d_{\ma L}$.
These points $\lbrace \ve x_j \rbrace_{j=1}^M$ have a special interpretation:
they are an $M$-point \emph{minimax optimal design}  on $\set D$:
\begin{equation}\label{eq:minimax}
	\set X(\set D, M, \ma L) \coloneqq \argmin_{\lbrace \ve x_j \rbrace_{j=1}^M\subset \set D}
	\max_{\ve x\in \set D}
	\min_{j = 1,\ldots, M}
	\|\ma L(\ve x - \ve x_j)\|_2;
\end{equation}
see proof in~\cite[Thm.~4.7]{CB05}.
By replacing the Lipschitz constant by the Lipschitz matrix
fewer points are required cover the domain
as illustrated in \cref{fig:cover}.
The same can be seen in bounds for the covering number.
\newsavebox{\Lmat}\savebox{\Lmat}{$\ma L=\begin{bsmallmatrix} \phantom{-}\sfrac{1}{4} & 0 \\-\sfrac{1}{4} & 1 \end{bsmallmatrix}$}
\DTLsettabseparator		\begin{filecontents*}{data/fig_cover_scalar.dat}
x	y	
-0.6686677797529386	0.6652016236538927	
0.6702697231060203	-0.6636102180321116	
-0.66486088961459	-0.6689895453532373	
0.6632884226676458	0.6706000602703281	
3.800657536298213e-05	-0.0016107864379552583	
0.6680531461348678	0.006997729474146657	
-0.6680400164042659	-0.0038012604283910252	
-0.005365376932752778	0.6658014099284093	
0.005392666938765086	-0.6689884548997165	
\end{filecontents*}
\DTLloaddb[noheader=false]{coverscalar}{data/fig_cover_scalar.dat}
\begin{filecontents*}{data/fig_cover_matrix.dat}
x	y	
-4.247944558852548e-05	0.8333014764654687	
4.381572104867137e-05	-0.8333306720172696	
-1.0308952715545682e-06	-2.9833609880778404e-05	
\end{filecontents*}
\DTLloaddb[noheader=false]{covermatrix}{data/fig_cover_matrix.dat}
\DTLmaketabspace		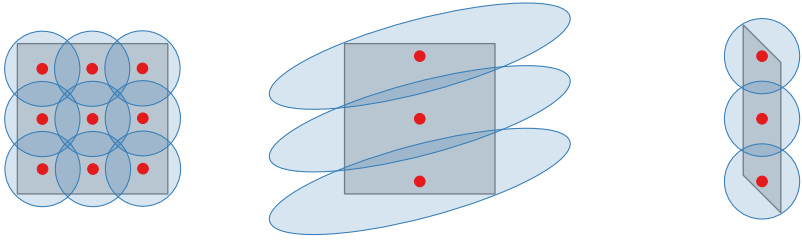
\begin{figure}
\centering
\begin{tikzpicture}
\begin{groupplot}[
		group style={group size = 3 by 1},
		width=0.5\linewidth,
		height = 0.5\linewidth,
		xmin = -2.05 , xmax = 2.05,
		ymin = -1.6 , ymax = 1.6,
		hide axis,
		x = 1cm,
		y = 1cm,
	]
	\nextgroupplot[title = {$L  = 1$, $\set D = [-1,1]^2$ 
		\vphantom{$\begin{bsmallmatrix} \sfrac{1}{4} \\ \sfrac{1}{4} \end{bsmallmatrix}$}
		},
		xmin = -1.3, xmax = 1.3,
		width = 0.3\linewidth,
	]
	
	\draw[line width=0.5pt, gray, fill=gray, fill opacity = 0.3] 
		(-1,-1) -- (-1,1) -- (1,1) -- (1,-1) -- cycle;

	\addplot[red, mark=*, only marks] table [x=x, y=y] {data/fig_cover_scalar.dat};
	\pgfplotsextra{\DTLforeach*{coverscalar}{\x=x, \y=y}{
		\filldraw[fill=colorbrewerA2,fill opacity=0.2, colorbrewerA2] (axis cs:\x,\y) circle [radius=0.5];
	}}

	\nextgroupplot[
		title = {$\ma L  = \begin{bsmallmatrix} \phantom{-}\sfrac{1}{4} & 0 \\ -\sfrac{1}{4} & 1\end{bsmallmatrix}$, $\set D = [-1,1]^2$},
		width=0.5\linewidth,
		xmin = -2.05 , xmax = 2.05,
		]
	
	\draw[line width=0.5pt, gray, fill=gray, fill opacity = 0.3] 
		(-1,-1) -- (-1,1) -- (1,1) -- (1,-1) -- cycle;

	\addplot[red, mark=*, only marks] table [x=x, y=y] {data/fig_cover_matrix.dat};
	\pgfplotsextra{\DTLforeach*{covermatrix}{\xx=x, \yy=y}{
		\draw[rotate around={14.872440898490666:(axis cs:\xx,\yy)}, colorbrewerA2, fill, fill opacity=0.2]
		(axis cs:\xx,\yy) ellipse [x radius = 2.06532429, y radius = 0.48418546];

	}}
	
	\nextgroupplot[
		title = {$\ma L  = \begin{bsmallmatrix} \phantom{-}\sfrac{1}{4} & 0 \\ -\sfrac{1}{4} & 1\end{bsmallmatrix}$, $\set D = \ma L [-1,1]^2$},
		xmin = -1.5, xmax =1.5,
		width = 0.35\linewidth,
	]

	\draw[line width=0.5pt, gray, fill=gray, fill opacity = 0.3] 
		(-0.25,-0.75) -- (-0.25,1.25) -- (0.25,0.75) -- (0.25,-1.25) -- cycle;

	\addplot[red, mark=*, only marks] table [x=x, y=y] {data/fig_cover_matrix.dat};
	\pgfplotsextra{\DTLforeach*{covermatrix}{\x=x, \y=y}{
		\filldraw[fill=colorbrewerA2,fill opacity=0.2, colorbrewerA2] (axis cs:\x,\y) circle [radius=0.5];
	}}

\end{groupplot}
\end{tikzpicture}

\caption{Using the Lipschitz matrix instead of the Lipschitz constant
reduces the number of $\epsilon$-balls required to cover a domain.
For an $\epsilon=0.5$ cover, 
9 points (dots) are required to cover a box (shaded quadrilateral)
with $\epsilon$-balls (shaded ellipses) using the Lipschitz constant (left).
In contrast, whereas only 3 points are required with the Lipschitz matrix (center, right).
The center plot shows the perspective of the Lipschitz matrix altering the metric for the space;
the right plot shows warping the domain to use the standard $\ell_2$ metric.
}
\label{fig:cover}
\end{figure}
 \begin{theorem}[{\cite[Thm.~14.2]{Wu17}}] \label{thm:cover}
	If $\set D\subset \R^m$ is convex
	and $\ma L$ is full rank
	\begin{equation}
		\left(\frac{1}{\epsilon}\right)^m
		\frac{\vol(\set D)}{\vol(\set B_1(\cdot; \ma L))}
		\le N_\epsilon(\set D; \ma L)
		\le 
		\left(\frac{3}{\epsilon}\right)^m
		\frac{\vol(\set D)}{\vol(\set B_1(\cdot; \ma L))},
	\end{equation}
	where $\vol$ denotes the Lebesgue measure in $\R^m$
	and $\set B_1(\cdot;\ma L)$ the unit ball in~$d_{\ma L}$.
\end{theorem}
This theorem implies $N_\epsilon(\set D;\ma L) = \order(\epsilon^{-m})$
when $\ma L$ is full rank.
If $\ma L$ is low-rank
we can apply this theorem by interpreting the Lipschitz matrix as warping the domain
\begin{equation}
	\ma L \set D \coloneqq \lbrace \ma L \ve x : \ve x \in \set D \rbrace 
\end{equation}
and equipping this domain with the standard Euclidean metric $d(\ve x, \ve y) = \|\ve x- \ve y\|_2$.
Then covering number is
\begin{equation}
	N_\epsilon(\set D; \ma L) = N_\epsilon(\ma L\set D;\ma I).
\end{equation}
When $\ma L$ is rank-$r$,
then $\ma L \set D$ is an $r$-dimensional subset of $\R^m$.
Applying \cref{thm:cover} on this $r$-dimensional problem, 
note $N_\epsilon(\ma L\set D;\ma I) = \order(\epsilon^{-r})$.

In the remainder of this section
we first establish a definition of uncertainty for Lipschitz functions
and then use this definition to derive complexity results for
approximation, integration, and optimization.
Then we discuss how the Lipschitz matrix yields
both asymptotic and non-asymptotic reductions to the complexity of
these tasks when compared to the Lipschitz constant.

\subsection{Uncertainty\label{sec:ibc:uncertainty}}
One use of the Lipschitz matrix is to provide constraints
on what values a function can take away from points where know its values.
We denote these \emph{point queries} consisting of a point $\ve x_j \in \set D$
and a response $y_j \in \R$ by $\set P_M = \lbrace \ve x_j , y_j \rbrace_{j=1}^M$.
This allows us to define the space of all Lipschitz functions with the Lipschitz matrix $\ma L$
that interpolate these point queries:
\begin{equation}
	\set L(\set D, \ma L, \set P_M) \coloneqq 
		\lbrace f \in \set L(\set D, \ma L): f(\ve x_j) = y_j \  \lbrace \ve x_j, y_j\rbrace \in \set P_M\rbrace.
\end{equation}
Using this notation we can then define the \emph{uncertainty set} at a point $\ve x\in \set D$
\begin{equation}\label{eq:uncertainty_set}
	\set U(\ve x; \ma L, \set P_M) \coloneqq
		\lbrace f(\ve x): f \in \set L(\set D, \ma L, \set P_M) \rbrace.
\end{equation}
This set is actually an interval for each $\ve x$.
From the definition of Lipschitz matrix continuity~\cref{eq:matrix_lipschitz},
for any $f \in \set L(\set D, \ma L, \set P_M)$ we have bounds
\begin{equation}\label{eq:lip_ineq}
	y_j - \|\ma L(\ve x - \ve x_j) \|_2 \le f(\ve x) \le y_j + \|\ma L(\ve x - \ve x_j)\|_2.
\end{equation}
As these bounds apply for each $\lbrace \ve x_j \rbrace_{j=1}^M$
these provide lower and upper bounds for each point $\ve x$:
\begin{align}
	\fmin(\ve x; \ma L, \set P_M) &\coloneqq \max_{j=1,\ldots, M} y_j - \|\ma L(\ve x - \ve x_j)\|_2, \\ 
	\fmax(\ve x; \ma L, \set P_M) &\coloneqq \min_{j=1,\ldots, M} y_j + \|\ma L(\ve x - \ve x_j)\|_2.
\end{align}
Since these are the only constraints on $f(\ve x)$,
the uncertainty set is an interval 
\begin{equation}\label{eq:uncertainty_interval}
	\set U(\ve x; \ma L, \set P_M) = \left[ \fmin(\ve x; \ma L, \set P_M), \ \fmax(\ve x; \ma L, \set P_M)\right]. 
\end{equation}

An important tool in our results is the \emph{central approximation}
\begin{equation}\label{eq:central_approx}
	\fmid (\ve x; \ma L, \set P_M) \coloneqq
		\frac12 \left( \fmin (\ve x) + \fmax (\ve x) \right).
\end{equation}
This function minimizes the worst-case pointwise error of all Lipschitz approximations
that interpolate the point set.
\begin{lemma}\label{lem:central}
	Given a domain $\set D \subset \R^m$,
	Lipschitz matrix $\ma L\in \R^{m\times m}$,
	and point queries $\set P_M = \lbrace \ve x_j, y_j\rbrace_{j=1}^M$
	then for any fixed $\widetilde f : \set D \to \R$ 
	\begin{equation}
		\sup_{f \in \set L(\set D, \ma L, \set P_M)}
			|f(\ve x) - \fmid(\ve x;\ma L, \set P_M)| \le
		\sup_{f \in \set L(\set D, \ma L, \set P_M)}
			|f(\ve x) - \widetilde f(\ve x)|.
	\end{equation}
\end{lemma}
\begin{proof}
	This follows immediately from the uncertainty interval~\cref{eq:uncertainty_interval}.
	Any other choice for $\widetilde f(\ve x)$
	would have an error greater than $\frac12[\fmax(\ve x) - \fmin(\ve x)]$ 
	for some choice of $f$;
	cf.,~\cite[pp.~12-13]{TW98}. 
\end{proof}

\subsection{Complexity Results\label{sec:ibc:complexity}}
The following results show 
the number of function queries are necessary to perform
approximation, integration, and optimization
to within a tolerance $\epsilon>0$ in the worst case
depends on the $\epsilon$-covering of the domain.

\subsubsection{Approximation}
The following theorem shows that by querying a function at a minimax optimal design,
the resulting central approximation yields the best approximation in the worst case.

\begin{theorem}\label{thm:approx}
	Given a domain $\set D \subset \R^m$ and a Lipschitz matrix $\ma L\in \R^{m \times m}$,
	the minimum number of point queries $M$ 
	such that 
	\begin{equation}\label{eq:approx_bound}
		\sup_{f \in \set L(\set D, \ma L)} 
			\inf_{\widetilde f: \set D \to \R}
			\max_{\ve x\in\set D} | f(\ve x) - \widetilde{f}(\ve x) | \le  \epsilon
		\quad \text{where} \quad
		\set P_M = \lbrace \ve x_j, f(\ve x_j)\rbrace_{j=1}^M	
	\end{equation}
	is the $\epsilon$ internal covering number, $M = N_\epsilon(\set D; \ma L)$.
	The optimal point queries correspond to the $M$-point minimax optimal design, $\set X(\set D, M, \ma L)$
	and the optimal approximation $\widetilde f$ is the central approximation $\fmid$. 
\end{theorem}
\begin{proof}
	Using \cref{lem:central} we can remove the middle optimization in~\cref{eq:approx_bound} 
	by replacing $\widetilde f$ with the central approximation $\fmid$:
	\begin{equation}
		\sup_{f \in \set L(\set D, \ma L)} 
			\inf_{\widetilde f: \set D \to \R}
			\max_{\ve x\in\set D} | f(\ve x) - \widetilde{f}(\ve x) | 
		=
			\sup_{f \in \set L(\set D, \ma L)} 
			\max_{\ve x\in\set D} | f(\ve x) - \fmid(\ve x; \ma L, \set P_M) |. 
	\end{equation}
	The greatest uncertainty comes when  $f(\ve x_j) = 0$ for $j=1,\ldots, M$
	and hence
	\begin{equation}
		\sup_{f \in \set L(\set D, \ma L)} 
			\max_{\ve x\in\set D} | f(\ve x) - \fmid(\ve x; \ma L, \set P_M) |
			= 			\max_{\ve x\in \set D} \min_{j=1,\ldots, M} \|\ma L(\ve x - \ve x_j)\|_2.
	\end{equation} 
	Then the minimum $M$ to make the right hand side above less than $\epsilon$
	is $N_\epsilon(\set D, \ma L)$
	and the optimal queries correspond $N_\epsilon(\set D, \ma L)$-point minimax optimal design.
\end{proof}

\subsubsection{Optimization}
The argument for the complexity of optimization is essentially that of approximation:
unless we can approximate $f$ within $\epsilon$,
we cannot globally optimize $f$ within $\epsilon$.

\begin{theorem}\label{thm:opt}
	Suppose $\set D\subset \R^m$ is compact and $\ma L\in \R^{m\times m}$.
	In the worst case, 
	the minimum number of samples to find the maximum $f\in \set L(\set D, \ma L)$
	to within $\epsilon$ is the $\epsilon$-internal covering number $N_\epsilon(\ma L\set D)$.
\end{theorem}
\begin{proof}
	As $\set D$ is compact and $f$ is Lipschitz,
	$f$ must have a finite maximizer $f^\star$
	where $f^\star = f(\ve x^\star)$ for at least one $\ve x^\star \in \set D$.
	We first establish an upper bound on the number of point queries required.
	Using $M = N_\epsilon(\set D, \ma L)$ point queries $\set P_M$
	in a minimax optimal design 
	there is a $\ve x_j$ in $\set P_M$ such that $\|\ma L(\ve x^\star - \ve x_j)\|_2 \le \epsilon$.
	By~\cref{eq:lip_ineq} 
	\begin{equation*}
		f(\ve x^\star) \le f(\ve x_j) + \|\ma L(\ve x - \ve x_j)\|_2 \le f(\ve x_j) + \epsilon
	\end{equation*}
	and hence $f(\ve x_j) \ge f^\star - \epsilon$.
	To show this upper bound is obtained,
	consider the function $f \in \set L(\set D, \ma L)$ 
	where $f(\ve x)\ge 0$ for all $\ve x\in \set D$,
	$f$ has optimizer $f(\ve x^\star) = 2\epsilon = f^\star$,
	and $f$ has minimum integrand.
	If we choose $M < N_\epsilon(\set D, \ma L)$ point queries
	we can adversarially choose $\ve x^\star$ such that $\|\ma L(\ve x - \ve x^\star)\|_2 > \epsilon$
	and consequently $f(\ve x_j)< \epsilon$ for all $j=1,\ldots, M$.
\end{proof}

\subsubsection{Integration}
This result parallels the one-dimensional Lipschitz case
presented by Traub and Werschulz~\cite[chap.~2]{TW98}.

\begin{theorem}\label{thm:integration}
	Suppose $\set D \subset \R^m$ be compact and $\ma L\in \R^{m\times m}$.
	Let $\phi^\star$ denote the quadrature rule from integrating the central approximation
	on an $M$-point minimax optimal design
	\begin{equation}\label{eq:opt_quadrature}
		\begin{split}
		&\phi^\star(f; M, \ma L) 
			\coloneqq \int_{\ve x\in \set D}
				\fmid(\ve x; \ma L, \set P_M) \D \ve x \\
		&\text{where} \quad
				\set P_M = \lbrace \ve x_j, f(\ve x_j)\rbrace_{j=1}^M
				\quad \text{and}
				\quad \ve x_j \in \set X(\set D, M, \ma L).
		\end{split}
	\end{equation}
	Let $\phi$ denote any other $M$ point quadrature rule integrating 
	an approximation $\widetilde f$,
	then 
	\begin{equation}\label{eq:int_bound}
		\sup_{f \in \set L(\set D, \ma L)} 
			\left| 
				\phi^\star(f) - \int_{\ve x \in \set D} f(\ve x) \D \ve x
			\right|
		\le \sup_{f\in \set L(\set D, \ma L)}
			\left| 
				\phi(f) - \int_{\ve x \in \set D} f(\ve x) \D \ve x
			\right|.
	\end{equation}
\end{theorem}
\begin{proof}
	From \cref{lem:central},
	the central approximation has the smallest pointwise worst case error
	of any approximation interpolating the point queries $\set P_M$.
	Thus
	\begin{equation}
		\sup_{f\in \set L(\set D, \ma L, \set P_M)}
			\left| \int_{\ve x\in \set D} \fmid(\ve x; \ma L, \set P_M) - f(\ve x) \D \ve x\right| 
		\le 
		\sup_{f\in \set L(\set D, \ma L, \set P_M)}
			\left| \int_{\ve x\in \set D} \widetilde f(\ve x) - f(\ve x) \D \ve x\right|.
	\end{equation}
	Then as the points $\set X(\set D, M, \ma L)$
	minimize the worse case error in the central approximation,
	$\phi^\star$ is the worst-case optimal quadrature rule.	
\end{proof}

\begin{corollary}
	In the setting of~\cref{thm:integration}, 
	\begin{equation}
		\left| 
			\phi^\star(f; N_\epsilon(\set D, \ma L), \ma L) 
			- \int_{\ve x\in \set D} f(\ve x) \D \ve x
		\right|
		\le \epsilon \cdot \vol(\set D).
	\end{equation}	
\end{corollary}
\begin{proof}
	Let $M = N_\epsilon(\set D, \ma L)$.
	From the definition of $\phi^\star$~\cref{eq:opt_quadrature}, 
	\begin{equation}
		\left| 
			\phi^\star(f; M, \ma L) - 
			\int_{\ve x\in \set D} f(\ve x) \D \ve x
		\right|
		\le 
		\int_{\ve x \in \set D}
			| \fmid(\ve x; \ma L, \set P_M) - f(\ve x)| \D \ve x.
	\end{equation}
	By \cref{thm:approx},
	$|\fmid(\ve x) - f(\ve x)| \le \epsilon$;
	integrating this bound yields the result.
\end{proof}

\subsection{Asymptotic Improvements\label{sec:ibc:asymptotic}}
The Lipschitz matrix provides two asymptotic improvements
over the Lipschitz constant.
If the Lipschitz matrix is rank-$r$,
then complexity of approximation, integration, and optimization 
grows like $\order(\epsilon^{-r})$---not $\order(\epsilon^{-m})$.
This also applies to ridge functions with ridge dimension $r$
by \cref{thm:ridge}.
Even if the Lipschitz matrix is full rank,
the Lipschitz matrix can still substantially reduce the covering number
compared to the Lipschitz constant.
Interpreting the Lipschitz matrix as transforming the domain,
from \cref{thm:cover} we have the bound 
\begin{equation}
	N_\epsilon(\set D; \ma L) = N_\epsilon(\ma L\set D; \ma I)  \le
		\frac{3^m}{\epsilon^{m}\vol(\set B_1(\cdot; \ma I)) } \vol(\ma L \set D).
\end{equation}
For the Lipschitz matrix and Lipschitz constant 
\begin{equation}
	\vol(\ma L\set D) = |\det \ma L| \cdot \vol(\set D)
	\quad \text{and} \quad
	\vol(L\set D) = L^m \cdot \vol(\set D).
\end{equation}
Hence if $|\det \ma L| \ll L^m$,
we have substantially reduced the covering number.
This is not uncommon:
for the test problems shown in \cref{tab:scaling},
$|\det \ma L|$ is multiple orders of magnitude smaller than $L^m$.

\begin{table}
\centering
\caption{Comparison of Lipschitz matrix, Lipschitz constant, and corresponding volumes
	for a variety of test functions~\cite{SB13} posed on the normalized domain $\set D = [-1,1]^m$
	with normalized outputs in $[0,1]$.
	Here we estimate $\ma L$ by minimizing the Frobenius norm 
	using $1024$ gradient samples
	including the $2^m$ corners and otherwise sampled randomly.
}
\label{tab:scaling}
\setlength\tabcolsep{4pt}
\begin{filecontents*}{data/tab_scaling.dat}
dimension	scalar_volume	matrix_volume	sing_max	sing_min	lam_max	lam_min	L_scalar	
6.0	0.0427248544045	8.38652306652e-09	0.304755924186	0.000725814426076	0.135549700576	0.000754177837822	0.295631088821	
6.0	1.6426663399	1.43480476112e-05	0.578507865968	0.00144747649851	0.3380838752	0.0014480209175	0.543118876577	
7.0	132.938872653	1.96202194923e-06	1.02219698971	0.00257677996265	0.785686498645	0.00272089031297	1.00542310272	
8.0	105.547330283	4.97142686205e-12	0.921243394238	4.61557642738e-06	0.497443881068	6.52741077821e-06	0.895160501512	
10.0	0.0585521452029	2.11323525617e-10	0.411478973615	0.00428607596521	0.201718646647	0.004299659252	0.376464722136	
\end{filecontents*}
\pgfplotstabletypeset[
	every head row/.style={before row=\toprule,after row=\midrule},
    every last row/.style={after row=\bottomrule},
	create on use/name/.style={
        create col/set list={
			\hspace*{-2pt}Golinski volume~\cite{Gol70},
			OTL circult~\cite{BS07a}, 
			piston~\cite{KZ98},
			borehole~\cite{HG83}, 
			wing weight~\cite{FSK08},
		},
    },
	create on use/ref/.style={
        create col/set list={
			\cite{Gol70},
			\cite{BS07a}, 
			\cite{KZ98},
			\cite{HG83}, 
			\cite{FSK08}
		},
    },
	columns/ref/.style = {column name = ref., string type, column type = {c},},
	columns/name/.style = {column name = test problem, string type, column type = {r}, },
	columns/dimension/.style = {int detect, column name = {\hspace*{-2pt}dim.\hspace*{-2pt}}  },
	columns/scalar_volume/.style = {column name = $\vol(L \set D)$, sci, precision = 1, sci zerofill, dec sep align},
	columns/matrix_volume/.style = {column name = $\vol(\ma L\set D)$, sci, precision = 1, sci zerofill, dec sep align},
	columns/sing_max/.style = {column name = $\sigma_\text{max}(\ma L)$, sci, precision = 1, sci zerofill, dec sep align},
	columns/sing_min/.style = {column name = $\sigma_\text{min}(\ma L)$, sci, precision = 1, sci zerofill, dec sep align},
	columns/lam_max/.style = {column name = $\lam_\text{max}(\ma L)$, sci, precision = 1, sci zerofill, dec sep align},
	columns/lam_min/.style = {column name = $\lam_\text{min}(\ma L)$, sci, precision = 1, sci zerofill, dec sep align},
	columns/L_scalar/.style = {column name = $L$, sci, precision = 1, sci zerofill, dec sep align},
	columns = {name, dimension, L_scalar, sing_max, sing_min, scalar_volume, matrix_volume},
]{data/tab_scaling.dat}
\end{table}
 
\subsection{Non-asymptotic Improvement\label{sec:ibc:nonasymptotic}}
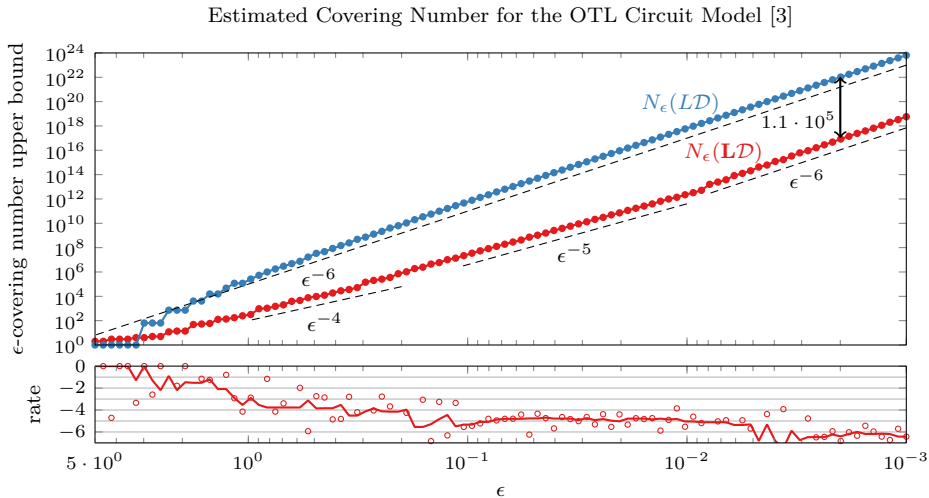
\begin{figure}
\centering 
\begin{tikzpicture}
\begin{groupplot}[
	group style = {
		group size = 1 by 2,
		vertical sep = 8pt,
	},
	xmode = log,
	ymode = log,
	width = 0.95\linewidth,
	height = 0.42\linewidth,
	ylabel = $\epsilon$-covering number upper bound,
	x dir = reverse,
	xlabel = $\epsilon$,
	ymin = 1,
	ymax = 1e24,
	xmin = 1e-3,
	xmax = 5e0,
	ytickten = {0,2,4,...,30},
	clip mode=individual,
	]
	
	\nextgroupplot[title = Estimated Covering Number for the OTL Circuit Model~\cite{BS07a},
		xticklabels = {,,},
		xlabel = {},
		]

	\addplot[colorbrewerA1, thick, mark=*, mark size =1pt] table [x=eps, y=N] {data/fig_covering_mat.dat}
		node[pos = 0.8, anchor = south east, yshift = -4pt] {$N_\epsilon(\ma L \set D)$} ;
	
	\addplot[black, densely dashed, domain = 2e-1:1e0] {1e2*(1/x)^4} 
		node [pos = 0.5, anchor = north] {$\epsilon^{-4}$};
	
	\addplot[black, densely dashed, domain = 1e-2:1.05e-1] {4e1*(1/x)^5} 
		node [pos = 0.5, anchor = north] {$\epsilon^{-5}$};
	
	\addplot[black, densely dashed, domain = 1e-3:8e-3] {7e-1*(1/x)^6} 
		node [pos = 0.5, anchor = north] {$\epsilon^{-6}$};
	
	\addplot[colorbrewerA2, thick, mark=*, mark size = 1pt] table [x=eps, y=N] {data/fig_covering_con.dat}
		node[pos = 0.8, anchor = south east, yshift = -4pt] {$N_\epsilon(L \set D)$} ;

	\addplot[black, densely dashed, domain = 1e-3:10] {1e5*(1/x)^6}
		 node [pos = 0.7, anchor = north west, yshift = 4pt] {$\epsilon^{-6}$};
	
	\addplot[black,<->, thick] coordinates {(2e-3,1e17) (2e-3,1e17*1.1e5)}
		node [pos = 0.3, anchor = east, xshift=2pt] {\scriptsize $1.1\cdot 10^{5}$};

	\nextgroupplot[
		height=0.2\linewidth,
		ylabel = rate, 
		ymode = linear,
		ymin = -7,
		ymax = 0,
		ytick = {-8,-6,...,-2,0},
		extra x ticks = {5e0},
		extra x tick labels = {$5\cdot 10^0$},
		]
	\addplot[black, opacity = 0.3] coordinates {(10,-1) (1e-4,-1)}; 
	\addplot[black, opacity = 0.3] coordinates {(10,-2) (1e-4,-2)}; 
	\addplot[black, opacity = 0.3] coordinates {(10,-3) (1e-4,-3)}; 
	\addplot[black, opacity = 0.3] coordinates {(10,-4) (1e-4,-4)}; 
	\addplot[black, opacity = 0.3] coordinates {(10,-5) (1e-4,-5)}; 
	\addplot[black, opacity = 0.3] coordinates {(10,-6) (1e-4,-6)}; 
	\addplot[black, opacity = 0.3] coordinates {(10,-7) (1e-4,-7)};

	\addplot[colorbrewerA1, mark=o, only marks, mark size = 1pt]
		 table [x = eps, y= slope] {data/fig_covering_mat_slope.dat};
	\addplot[colorbrewerA1, thick] 
		table [x = eps, y= median] {data/fig_covering_mat_slope.dat};

\end{groupplot}
\end{tikzpicture}

\caption{An upper bound on the covering number of the transformed domain.
	Here we estimate the covering number by counting the number of $\epsilon$-balls
	whose centers are on a grid with spacing $2\epsilon/\sqrt{m}$
	that intersect the transformed domain;
	when the number of grid points exceeds $10^4$, 
	random samples are used to estimate the number of $\epsilon$-balls intersecting the domain.
	The bottom plot shows the estimated growth rate of $N_\epsilon(\ma L\set D)$
	estimated using a finite difference (dots) and a 7-point median smoothed rate (line). 
	The asymptotic separation of $1.1\cdot 10^5$ matches the ratio $\vol(L\set D)/\vol(\ma L\set D)$.
}
\label{fig:covering}
\end{figure}

Before the asymptotic limit,
the Lipschitz matrix can slow the growth of covering number.
Denoting the singular values of $\ma L$ in non-increasing order as $\sigma_j(\ma L)$,
when $\epsilon\gg \sigma_j(\ma L)$ then $\ma L\set D$ is effectively $j$-dimensional
as a single $\epsilon$-ball can cover the dimensions $j+1,\ldots,m$.
This temporally slows the growth of the covering number
as illustrated in \cref{fig:covering}. 
  \section{Computing a Lipschitz Matrix\label{sec:estimate}}
For a particular function $f$, 
our goal is to identify the smallest Lipschitz matrix $\ma L$
such that $f\in \set L(\set D, \ma L)$.
As discussed in \cref{sec:intro:smallest},
we must choose an ordering of positive semidefinite matrices.
Here we minimize the Frobenius norm of the Lipschitz matrix
yielding the program
\begin{equation}\label{eq:lip_all}
	\begin{split}
	\min_{\ma L\in \mathbb{S}_+^{m\times m}} & \ \|\ma L\|_\fro^2  \\
	\text{such that} & \
	|f(\ve x_i) - f(\ve x_j)| \le \|\ma L(\ve x_i - \ve x_j)\|_2 
	\quad \forall \ve x_i, \ve x_j \in \set D.
	\end{split}
\end{equation}
Ideally given $f$, 
we would identify this Frobenius-norm minimizing Lipschitz matrix in closed form.
For most functions, this is infeasible.
Instead, we approximate this program
by discretizing the constraints 
using a finite number of point and/or gradient queries.
This yields a semidefinite program to identify the Lipschitz matrix.
Although tempting to enforce a low-rank constraint on the Lipschitz matrix in this program,
we show this can yield in non-informative results.
For the special case of quadratic functions, 
we provide a finite program exactly solving the full program~\cref{eq:lip_all}.
We finally illustrate how this discretized program
converges with increasing queries.

\subsection{Semidefinite Program for a Lipschitz Matrix}
Both point queries 
$\set P_M = \lbrace \ve x_j , f(\ve x_j)\rbrace_{j=1}^M$
and gradient queries $\set G_N = \lbrace \nabla f(\ve x_k') \rbrace_{k=1}^N$
provide constraints on the minimum Frobenius norm Lipschitz matrix.
Recalling from~\cref{eq:grad_order}
that the outer product of gradients is bounded above by 
$\ma L^\trans \ma L$, 
the discretized minimum Frobenius-norm Lipschitz matrix solves
\begin{equation}\label{eq:Lopt_finite}
	\begin{aligned}
	\min_{\ma L\in \mathbb{S}_+^{m\times m}} && \ \|\ma L\|_\fro^2   \\
	\text{such that} && 
	|y_i - y_j|^2 &\le \|\ma L(\ve x_i - \ve x_j)\|_2^2
	\quad & \lbrace \ve x_i,y_i \rbrace, \lbrace \ve x_j, y_j \rbrace &\in \set P_M;
	\\
	&& \ve g_k \ve g_k^\trans &\preceq \ma L^\trans \ma L
	& \ve g_k &\in \set G_N.
	\end{aligned}
\end{equation}
This formulation has difficult non-convex quadratic constraints.
Instead, we reformulate~\cref{eq:Lopt_finite}
in terms of the squared Lipschitz matrix 
$\ma H \coloneqq \ma L^\trans \ma L$ and instead solve
\begin{equation}\label{eq:Hopt_finite}
	\begin{aligned}
	\min_{\ma H\in \mathbb{S}_+^{m\times m}} && \ \Trace \ma H   \\
	\text{such that} && 
	|y_i - y_j|^2 &\le (\ve x_i - \ve x_j)^\trans \ma H(\ve x_i - \ve x_j) 
	& \lbrace \ve x_i,y_i \rbrace, \lbrace \ve x_j, y_j \rbrace &\in \set P_M;
	\\
	&& \ve g_k \ve g_k^\trans &\preceq \ma H
	& \ve g_k &\in \set G_N.
	\end{aligned}
\end{equation}
This, unlike~\cref{eq:Lopt_finite}, is convex semidefinite program.
After parametrizing the space of symmetric matrices,
this program has three sets of constraints:
linear inequality constraints from the point queries,
semidefinite constraints from the gradient queries,
and the semidefinite constraint $\ma 0 \preceq \ma H$. 
Our numerical experiments
use CVXOPT~\cite{cvxopt} to solve~\cref{eq:Hopt_finite}
and take $\ma L$ to be the symmetric positive definite square root of $\ma H$.

\subsection{Semidefinite Program for an $\epsilon$-Lipschitz Matrix}
We use a similar approach working in terms of the squared Lipschitz matrix $\ma H = \ma L^\trans \ma L$
to compute an $\epsilon$-Lipschitz matrix.
However, because $\epsilon$-Lipschitz functions are not differentiable, 
gradient queries $\set G_N$ do not constrain this matrix.
This leaves only point queries which by~\cref{eq:Lmat_epsilon} must satisfy 
\begin{equation}
	|f(\ve x_1) - f(\ve x_2)| \le \epsilon + \|\ma L(\ve x_1 - \ve x_2)\|_2.
\end{equation}
Replacing the point query constraint in~\cref{eq:Hopt_finite}
with this expression yields 
\begin{equation}\label{eq:Hopt_epsilon}
	\begin{aligned}
	\min_{\ma H\in \mathbb{S}_+^{m\times m}} & \ \Trace \ma H   \\
	\text{such that} & \  
	\max(|y_i \!-\! y_j| \!- \! \epsilon,0)^2 \le \! (\ve x_i - \ve x_j)^{\!\trans} \ma H(\ve x_i - \ve x_j) 
	\ \  \lbrace \ve x_i,y_i \rbrace ,\lbrace \ve x_j, y_j\rbrace \! \in\! \set P_M .\!\!\!
	\end{aligned}
\end{equation}

\subsection{Low-rank Solutions}
It is tempting to impose a rank constraint on $\ma H$, 
and consequently $\ma L$,
to reduce complexity.
Unfortunately, a low-rank constraint can yield misleading results.
For example, suppose we have point queries $\lbrace \ve x_j, y_j\rbrace_{j=1}^M = \set P_M$
where $\lbrace \ve x_j\rbrace_{j=1}^M$ is in general position.
Then for almost every vector $\ve a \in \R^m$,
the projections onto $\Span \ve a$ are distinct: $\ve a^\trans \ve x_i \neq \ve a^\trans \ve x_j$ for $i \ne j$.
Thus, for almost every $\ve a$ there is a rank-1 Lipschitz matrix
whose range is $\Span \ve a$:
\begin{equation}
	\ma L_{\ve a} = \begin{bmatrix}
			\ma 0 \\
			\alpha \ve a^\trans
		\end{bmatrix}\in \R^{m \times m},
	\quad 
	\alpha = \max_{i \ne j} \frac{ |y_i-y_j|}{|\ve a^\trans \ve x_i - \ve a^\trans \ve x_j|}.
\end{equation}
Hence regardless of the actual structure of $f$
we have likely mistakenly identified it as a one-dimensional ridge function.
Gradient constraints for the Lipschitz matrix do not fail in a similar way.
If $\Range\lbrace \nabla f(\ve x_k') \rbrace_{j=1}^N$ is $n$-dimensional,
then $\ma H$, and consequently $\ma L$, must be at least rank-$n$.

\subsection{Using the Determinant\label{sec:estimating:det}}
In \cref{sec:ibc:asymptotic} we saw that the complexity of approximation, integration, and optimization
are proportional to $|\det \ma L|$.
Why not use the determinant as the objective function instead 
of the Frobenius norm in~\cref{eq:lip_all}?
There are two important reasons.
We loose convexity: $\|\ma L\|_\fro^2 = \Trace \ma H$ is convex 
whereas $|\det \ma L|^2 = \det \ma H$ is concave.
The other is that the determinant yields uninformative Lipschitz matrices.
As illustrated in the previous subsection, 
given finite point queries we can always find a low-rank Lipschitz matrix.
As the determinant of this matrix is zero
it is an optimal, but uninformative, Lipschitz matrix.

\subsection{Quadratic Functions\label{sec:estimating:quadratic}}
In the case of a quadratic function
we can show that there are only finite number of active gradient constraints
of the continuous problem~\cref{eq:lip_all}
yielding a finite semidefinite program for its Lipschitz matrix. 
Suppose $f$ is a quadratic function:
\begin{equation}\label{eq:quadratic}
	\set D = [-1,1]^m \quad
	f:\set D \to \R 
	\quad \text{where} \quad
	f(\ve x) = \ve x^\trans \ma A \ve x + \ve b^\trans \ve x,
\end{equation}
$\ma A\in \R^{m\times m}$, and $\ve b\in \R^m$
whose gradient is $\nabla f(\ve x) = (\ma A + \ma A^{\!\trans}) \ve x + \ve b$.
The following theorem bounds the gradient outer-product
above by points on the corners of the domain:
those points where each coordinate takes on the value $-1$ or $1$.

\begin{theorem}\label{thm:corner}
	Suppose $f:[-1,1]^m\to \R$ is a quadratic function
	as in~\cref{eq:quadratic}.
	For any $\ve x \in \set D$ there is a point $\ve c$ on a corner of $\set D$ such that
	\begin{equation}
		\nabla f(\ve x) \nabla f(\ve x)^\trans \preceq 
		\nabla f(\ve c) \nabla f(\ve c)^\trans.
	\end{equation}
\end{theorem}
\begin{proof}
	Denote $\ma A_S = \ma A + \ma A^\trans$.
	Consider the difference of these two gradient outer-products
	for any $\ve c$ on the corner of $\set D$ and $\ve x\in \set D$:
	\begin{multline}
		\nabla f(\ve c) \nabla f(\ve c)^\trans - \nabla f(\ve x) \nabla f(\ve x)
		=\\ \ma A_S (\ve c - \ve x)\ve b^\trans + \ve b (\ve c - \ve x)^\trans \ma A_S
			+ \ma A_S( \ve c \ve c^\trans - \ve x\ve x^\trans) \ma A_S.
	\end{multline}
	As $\ve c$ is on the corner of $[-1,1]^m$, then $\ve c \ve c^\trans \succeq \ve x \ve x^\trans$ and	
	\begin{equation}
		\nabla f(\ve c) \nabla f(\ve c)^\trans - \nabla f(\ve x) \nabla f(\ve x)
		\succeq  \ma A_S (\ve c - \ve x)\ve b^\trans + \ve b (\ve c - \ve x)^\trans \ma A_S.
	\end{equation}
	As $\ve c$ is on a corner, 
	there is a $\ve c$ such that
	the entires of the entries of $\ve c - \ve x$
	can have any combination of signs $\lbrace +,0\rbrace$ or $\lbrace -, 0\rbrace$.
	Thus we can choose a corner such that	 
	$\ma A_S(\ve c - \ve x) \ve b^\trans$ has nonnegative entries.
	Then the right hand side above is nonnegative and we conclude
	\begin{equation}
		\nabla f(\ve c) \nabla f(\ve c)^\trans - \nabla f(\ve x) \nabla f(\ve x) \succeq \ma 0.
	\end{equation}
\end{proof}

As a result of this theorem, 
\emph{the} Frobenius-norm minimizing Lipschitz matrix for a quadratic function
is the solution to a finite-dimensional semidefinite program
\begin{equation}\label{eq:Hopt_quad}
	\begin{aligned}
	\argmin_{\ma H\in \mathbb{S}_+^{m\times m}} & \ \Trace \ma H   \\
	\text{such that}  
	& \ \ve g_k \ve g_k^\trans \preceq \ma H
	\qquad   \ve g_k \in \lbrace \nabla f(\ve c_k)\rbrace_{k=1}^{2^m}
	\end{aligned}
\end{equation}
where $\lbrace \ve c_k\rbrace_{k=1}^{2^m}$ are the $2^m$ corners of $\set D = [-1,1]^m$.

\subsection{Convergence in Queries}
With increasing point and gradient queries, 
the solution the finite-constraint semidefinite program~\cref{eq:Lopt_finite}
converges to the continuous-constraint program~\cref{eq:lip_all}
if the function $f$ is Lipschitz.
Let $\ma H_{M,N}$ denote the solution to \cref{eq:Hopt_finite} with $\set P_M$ and $\set G_N$.
If $\set P_M \subset \set P_{M+1}$ and $\set G_N\subset \set G_{N+1}$
then $\Trace \ma H_{M,N} \le \Trace \ma H_{M+1,N}$ and $\Trace \ma H_{M, N} \le \Trace \ma H_{M, N+1}$.
As $f$ is Lipschitz there exists a solution $\ma H$ to the continuous-constraint problem 
and $\Trace \ma H_{M,N} \le \Trace \ma H$.
Then as $\Trace \ma H_{M,N} \to\Trace\ma H$,
$\ma H_{M,N} \to \ma H$.
 
\subsection{Convergence Rate}
\begin{figure}
\pgfplotstableread{data/fig_convergence_samp.dat}\samp
\pgfplotstableread{data/fig_convergence_grad.dat}\grad
\pgfplotstableread{data/fig_convergence_time_samp.dat}\tsamp
\pgfplotstableread{data/fig_convergence_time_grad.dat}\tgrad
\pgfplotstableread{data/fig_convergence_time_svd.dat}\tsvd
\centering

\begin{tikzpicture}
\begin{groupplot}[
	group style = {group size = 1 by 2, 
		horizontal sep = 14pt,
		vertical sep = 10pt},
	width = 0.93\textwidth,	
	height = 0.35\textwidth,
	xmode = log,
	ymode = log,
	]

	\nextgroupplot[title = Comparison of Sample and Gradient Based Lipschitz Matrix Estimates,
		ymin = 5e-2,
		ymax = 1,
		xmin = 1,
		xmax = 1e4,
		xticklabels = {,,},
		ymode = log,
		ylabel = {$\| \ma H - \hma H \|_\fro/\|\hma H\|_\fro$},
	]
	
	\addplot[colorbrewerA2, thick] table [x=M, y = p50] {\samp}
		node [pos=0.8, anchor = south west, yshift =0pt,rotate=-9] {samples};

	\addplot[draw=none, name path = mismatch_samp_lb] table [x=M, y = p25] {\samp};
	\addplot[draw=none,name path = mismatch_samp_ub] table [x=M, y = p75] {\samp};
	\addplot[colorbrewerA2, opacity = 0.3] fill between [of=mismatch_samp_lb and mismatch_samp_ub];

	\addplot[colorbrewerA1, thick] table [x=N, y = p50] {\grad}
		node [pos=0.7, anchor = north east, yshift=2pt, rotate=-9] {gradients};
	\addplot[draw=none, name path = mismatch_grad_lb] table [x=N, y = p25] {\grad};
	\addplot[draw=none, name path = mismatch_grad_ub] table [x=N, y = p75] {\grad};
	\addplot[colorbrewerA1, opacity = 0.3] fill between [of=mismatch_grad_lb and mismatch_grad_ub];

	\addplot[colorbrewerA3, thick] table [x expr = \thisrow{N}*7, y = p50] {\grad}
		node [pos=0.5, anchor = south west, rotate=-9] {FD-gradients};
	\addplot[draw=none, name path = mismatch_fdgrad_lb] table [x expr=\thisrow{N}*7, y = p25] {\grad};
	\addplot[draw=none, name path = mismatch_fdgrad_ub] table [x expr=\thisrow{N}*7, y = p75] {\grad};
	\addplot[colorbrewerA3, opacity = 0.3] fill between [of=mismatch_fdgrad_lb and mismatch_fdgrad_ub];

	\addplot[black, densely dashed, domain = 1:1e4] {0.25*x^(-1/6)} 
		node [pos = 0.3, anchor = north, rotate=-9] {$\order(N^{-1/6})$};

	\nextgroupplot[
		xmin = 1,
		xmax = 1e4,
		ymin = 5e-3,
		ymax = 1e3,
		ylabel = wall time (seconds),
		xlabel = {quantity of data: $M$ samples or $N$ gradients},
		ytickten = {-2,...,3},
		]
	\addplot[colorbrewerA2, thick] table [x=M, y = p50] {\tsamp}
		node [pos=0.75, anchor = south east, rotate=25] {samples};
	\addplot[draw=none, name path = time_samp_lb] table [x=M, y = p25] {\tsamp};
	\addplot[draw=none, name path = time_samp_ub] table [x=M, y = p75] {\tsamp};
	\addplot[colorbrewerA2, opacity = 0.3] fill between [of=time_samp_lb and time_samp_ub];
	
	\addplot[black, densely dashed, domain = 2:1e3] {2e-4*x^2} 
		node [pos = 0.95, anchor = north, xshift=2pt, rotate=25] {$\order(M^2)$};

	\addplot[colorbrewerA1, thick] table [x=N, y = p50] {\tgrad}
		node [pos=1, anchor = south east, xshift=0pt, yshift=-3pt,rotate=12] {gradients};
	\addplot[draw=none, name path = time_grad_lb] table [x=N, y = p25] {\tgrad};
	\addplot[draw=none, name path = time_grad_ub] table [x=N, y = p75] {\tgrad};
	\addplot[colorbrewerA1, opacity = 0.3] fill between [of=time_grad_lb and time_grad_ub];
	
	\addplot[black, densely dashed, domain = 1:1e4] {2e-3*x} 
		node [pos = 0.7, anchor = north,rotate=12] {$\order(N)$};

\end{groupplot}
\end{tikzpicture}

\caption{The convergence of a discrete estimate
	of the Frobenius norm minimizing squared Lipschitz matrix $\ma H$
	to it true value $\hma H$
	for a least squares quadratic approximation to the 6-parameter OTL circuit function
	on a $8^6$ point grid.
	Finite difference gradients 
	are forward differences in each unit direction
	from one central point requiring seven function evaluations.
	For each number of queries $M$ and $N$
	we randomly sample the domain with uniform probability.
	One hundred repetitions are performed,
	with the median shown by a solid line and 
	the shaded region enclosing the $25$th to $75$th percentile.
	The top plot shows the error in the squared Lipschitz matrix
	while the bottom plot shows the wall clock time 
	on a Xeon E5-2620 clocked at 2.1 GHz using OpenBLAS with one thread
	and CVXOPT 1.2.5.
}
\label{fig:convergence}
\end{figure}
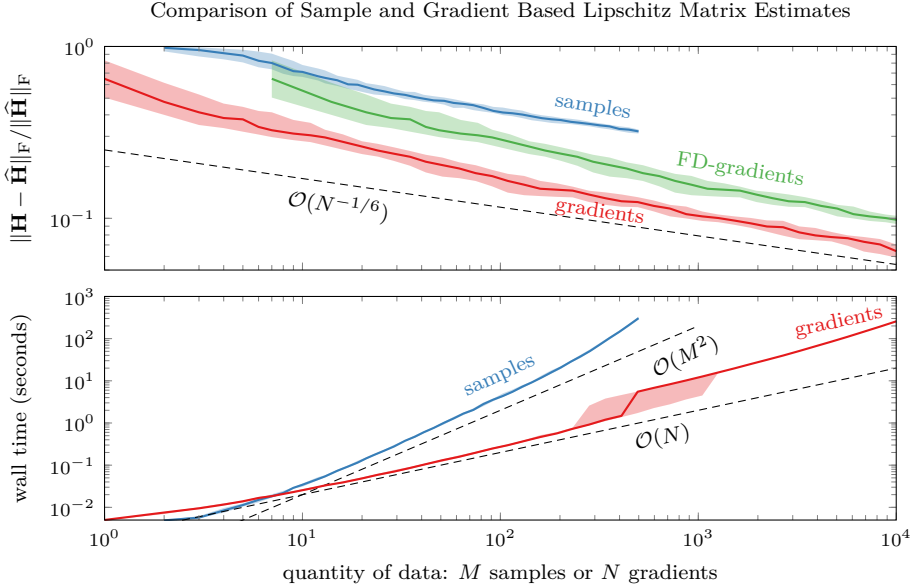

An important practical question is
how fast do finite query approximations converge 
to the continuous-constraint Lipschitz matrix~\cref{eq:lip_all}?
Unfortunately if we use random sampling this rate can be very slow.
\Cref{fig:convergence} shows a quadratic approximation to the OTL Circuit function.
As this is a quadratic function, we know the Lipschitz matrix 
is determined by the gradient at the corners.
The probability of querying within $\epsilon$ of these points
when sampling randomly is $\order(\epsilon^{-m})$;
hence convergence with gradient queries is $\order(\epsilon^{-m})$
and convergence with point queries is $\order(\epsilon^{-2m})$.
This is very slow!
However, this slow convergence
is not surprising as the Lipschitz matrix tracks the \emph{maximum} rate of change.

Although this result is distressing,
hope is not lost.
Many problems in science and engineering
are approximately monotonic over their domains.
Hence we can employ a similar argument to \cref{thm:corner}
to justify that querying gradients at the corners of the domain yields an 
accurate approximation of the Lipschitz matrix.
We employ this strategy throughout the remainder of this manuscript
to estimate the Lipschitz matrix, 
potentially in combination with other random point or gradient queries.

 \section{Ridge Approximation Error Bound\label{sec:ridge}}
In the introduction, \cref{thm:ridge} 
showed that functions with low-rank Lipschitz matrices ridge functions.
Here we provide a related result that bounds the error made when 
approximating Lipschitz function by a ridge function.
This result is analogous to the error bound associated with 
the MeGO matrix~\cite[Thm.~4.3]{Con15}.

\begin{theorem}\label{thm:error_bound}
	Suppose $f\in \set L(\set D, \ma L)$,
	$\widetilde{f} \in \set L(\set D, \ma L\ma U\ma U^\trans)$
	where $\ma U\in \R^{m\times n}$ and $\ma U^\trans \ma U = \ma I$,
	and $\lbrace \ve x_j \rbrace_{j=1}^M \subset \set D$, 
	then 
	\begin{multline}
			\max_{\ve x \in \set D} | f(\ve x) - \widetilde{f}(\ve x)| \le 
			\max_{j=1,\ldots, M} |f(\ve x_j ) - \widetilde f(\ve x_j)| + \\ 
			2\max_{\ve x \in \set D} \min_{j=1,\ldots, M} \| \ma L \ma U\ma U^\trans(\ve x_j - \ve x)\|_2
		 + \sigma_{\max}(\ma L(\ma I - \ma U\ma U^\trans)) \cdot \diam(\set D)
	\end{multline}
	where $\sigma_{\max}(\ma A)$ denotes the largest singular value of $\ma A$ and
	$\diam(\set D)$ denotes the diameter of the set $\set D$;
	$\diam (\set D) \coloneqq \max_{\ve x_1, \ve x_2 \in \set D} \|\ve x_1 - \ve x_2\|_2$.
\end{theorem}
\begin{proof}
Suppose $\ve x \in \set D$ is fixed.
Let $j = \argmin_k\|\ma L \ma U \ma U^\trans (\ve x_k - \ve x)\|_2$
and define  $\epsilon \coloneqq \max_{j=1,\ldots, M}|f(\ve x_j) - \widetilde{f}(\ve x_j)|$.
Inserting an additive identity,
\begin{align}
	|f(\ve x) \!-\! \widetilde{f}(\ve x)| &= 
			| f(\ve x) - \widetilde{f}(\ve x) 
				+ [f(\ve x_j) - \widetilde{f}(\ve x_j)] 
				- [f(\ve x_j) - \widetilde{f}(\ve x_j)] 
			  |  \\
		&\le  |f(\ve x)   -  f(\ve x_j)| + | \widetilde{f}(\ve x)  -  \widetilde{f}(\ve x_j)| + \epsilon.
\end{align}
Invoking each function's Lipschitz matrix and again inserting the identify,
\begin{align}
	|f(\ve x) \!-\! \widetilde{f}(\ve x)| 
		&\le \| \ma L(\ve x - \ve x_j)\|_2 + \| \ma L\ma U\ma U^\trans (\ve x - \ve x_j)\|_2  + \epsilon\\
		&\le \|\ma L( \ma U\ma U^\trans \! + \ma I - \ma U\ma U^\trans)(\ve x - \ve x_j)\|_2
			+ \|\ma L(\ma U\ma U^\trans)(\ve x - \ve x_j)\|_2 +\epsilon.
\end{align}
Finally, using the triangle inequality in the last term,
\begin{align}
	|f(\ve x) \!-\! \widetilde{f}(\ve x)| 
		&\le 2\|\ma L\ma U\ma U^\trans(\ve x - \ve x_j)\|_2 
			+ \|\ma L(\ma I - \ma U\ma U^\trans)(\ve x - \ve x_j)\|_2 + \epsilon. 
\end{align}
Then defining $\delta \coloneqq \max_{\ve x \in \set D} 
\min_{j=1,\ldots, M} \|\ma L\ma U\ma U^\trans(\ve x - \ve x_j)\|_2$
and bounding the second term by the diameter, we obtain a result independent of $\ve x$
\begin{align}
		|f(\ve x) - \widetilde f(\ve x)| & \le \epsilon + 2\delta + \sigma_{\max}(\ma L(\ma I - \ma U\ma U^\trans))\cdot\diam(\set D).
\end{align} 
\end{proof}

We can then use this theorem to 
motivate a particular choice for $\ma U$ 
and function evaluations $\ve x_j$.
First note that each term in this theorem has an important interpretation:
\begin{align*}
	&\max_{j} |f(\ve x_j) - \widetilde f(\ve x_j)| & 
	&\text{\emph{function approximation error} on $\lbrace \ve x_j \rbrace_{j=1}^M$, } \\
	&\max_{\ve x\in \set D} \min_{j} \|\ma L\ma U \ma U^\trans (\ve x_j - \ve x)\|_2 &
	&\text{\emph{dispersion} of $\lbrace \ve x_j \rbrace_{j=1}^M$ in $\set D$,} \\
	&\sigma_{\max}(\ma L(\ma I -\ma U\ma U^\trans)) &
	&\text{ \emph{Lipschitz matrix approximation error}}.
\end{align*}
To minimize the Lipschitz matrix approximation error
we should choose the 
leading eigenvectors of $\ma L$
as the columns of $\ma U$
(or right singular vectors of $\ma L$ is not symmetric positive definite).
With $\ma U$ selected,
we can minimize dispersion 
by constructing a minimax design of experiments 
$\set X(\set D, M, \ma L\ma U\ma U^\trans)$ as in \cref{eq:minimax}.
As discussed in \cref{sec:ibc:uncertainty},
the number of function evaluations required to obtain a particular dispersion (the covering number)
no longer grows with the dimension of the domain,
but instead the rank of the Lipschitz matrix $\ma L \ma U\ma U^\trans$.
If $f$ is a ridge function then we can make the function approximation error zero;
otherwise, this is not the case.
If two coordinates have the same projection,
$\ma U^\trans \ve x_i = \ma U^\trans \ve x_j$,
and  $f(\ve x_i) \ne f(\ve x_j)$,
then the function approximation error is at least
$|f(\ve x_i) - f(\ve x_j)|/2$.

Before concluding, 
we note two corollaries
when the function approximation error is zero.

\begin{corollary}\label{cor:error_bound:identity}
	In the setting of \cref{thm:error_bound}, let $\ma U$ be the identity matrix
	and choose $\widetilde f$ such that $f(\ve x_j) = \widetilde f(\ve x_j)$, 
	then 
	\begin{equation}
		\max_{\ve x\in \set D} |f(\ve x) - \widetilde f(\ve x)| 
		\le 2\max_{\ve x \in \set D} \min_{j = 1,\ldots, M} \|\ma L(\ve x_j - \ve x)\|_2.
	\end{equation}
\end{corollary}

\begin{corollary}\label{cor:error_bound:ridge}
	In the setting of \cref{thm:error_bound}, let $f$ be a ridge function
	and $\ma U$ be a basis for the range of $\ma L$.
	Taking $\widetilde f$ such that $f(\ve x_j) = \widetilde f(\ve x_j)$, 
	then 
	\begin{equation}
		\max_{\ve x\in \set D} |f(\ve x) - \widetilde f(\ve x)| 
		\le 2\max_{\ve x \in \set D} \min_{j = 1,\ldots, M} \|\ma L\ma U\ma U^\trans(\ve x_j - \ve x)\|_2.
	\end{equation}
\end{corollary}

 \section{Design of Experiments\label{sec:design}}
Given a particular function, 
where should we evaluate it to provide the most information?
This is the subject of the \emph{design of computer experiments}~\cite{SWN03},
a subfield of \emph{experimental design} (see, e.g.,~\cite{Fed72})
distinguished by the assumption that observations are deterministic;
e.g., $f(\ve x)$ returns only one value. 
From the ridge approximation error bound in~\cref{thm:error_bound}
and the earlier complexity results in \cref{sec:ibc:complexity},
a good experimental design of points $\lbrace \ve x_j \rbrace_{j=1}^M \subset \set D$
should minimize the \emph{dispersion} in the Lipschitz matrix metric:
\begin{equation}
	\delta(\lbrace \ve x_j \rbrace_{j=1}^M, \set D, \ma L)
		\coloneqq
		\max_{\ve x \in \set D} 
		\min_{j=1,\ldots, M} \|\ma L( \ve x - \ve x_j)\|_2.
\end{equation}
Minimizing the dispersion yields in a \emph{minimax optimal design}; cf.~\cref{eq:minimax}:
\begin{equation}\label{eq:minimax2}
	\set X(\set D, M, \ma L) = \argmin_{\lbrace \ve x_j \rbrace_{j=1}^M\subset \set D}
	\max_{\ve x\in \set D}
	\min_{j = 1,\ldots, M}
	\|\ma L(\ve x - \ve x_j)\|_2.
\end{equation}
There is a substantial body of literature on solving this problem;
see, e.g., \cite{Pro17} for a recent review.
Here we provide a motivating example 
illustrating why an minimax design in the Lipschitz metric
is important.
Then we provide a brief description of how we construct
locally optimal low-dispersion designs 
and illustrate the performance of this algorithm.

\subsection{A Motivating Example}
\begin{figure}
\centering
\begin{tikzpicture}
\begin{groupplot}[
		group style = {group size = 2 by 1, horizontal sep=10pt},
		height = 0.3\linewidth,
		width = 0.5 \linewidth,
	]

	\nextgroupplot[
		title = Prototypical Design,
		xlabel = {$\ve 1^\trans \ve x$},
		ymin = -0.5, 
		ymax = 2.5,
		ytick = {0,1,2},
		yticklabels = {minimax, LHS, random},
		xmin = -10.2,
		xmax = 10.2,
		xtick = {-10,-8,...,10},
		axis y line=left,
		y axis line style={opacity=0},
		axis x line=bottom,
		ytick style={draw=none},
		axis line style={-},
	]

		\addplot[blue, only marks, mark = |, thick, mark size = 5pt]
			table [x expr = \thisrow{x}*sqrt(10), y expr = 2] {data/fig_sine_rand.dat};
		\addplot[black, thick]
			coordinates { (-10,2) (10,2)};	
		\addplot[green, only marks, mark = |, thick, mark size = 5pt]
			table [x expr = \thisrow{x}*sqrt(10), y expr = 1] {data/fig_sine_lhs.dat};
		\addplot[black, thick]
			coordinates { (-10,1) (10,1)};	
		\addplot[red, only marks, mark = |, thick, mark size = 5pt]
			table [x expr = \thisrow{x}*sqrt(10), y expr = 0] {data/fig_sine_minimax.dat};
		\addplot[black, thick]
			coordinates { (-11,0) (10,0)};

	\nextgroupplot[
		ymin = -0.5,
		ymax = 2.5,
		ytick = {0,1,2},
		yticklabels = {,,},
		xmin = 0,
		xmax = 30,
		xtick = {0,5,10,...,30},
		xlabel = dispersion,
		y dir = reverse,
		title = Dispersion of $100$ Designs,
		axis y line=left,
		y axis line style={opacity=0},
		axis x line=bottom,
		ytick style={draw=none},
		axis line style={-},
	]

	\addplot[blue, only marks, mark=*, mark size = 0.5pt]
		table [x expr = 10*\thisrow{y}, y expr= 1.7*\thisrow{x}] {data/fig_sine_swarm_rand.dat};
	
	\addplot[green, only marks, mark=*, mark size = 0.5pt]
		table [x expr = 10*\thisrow{y}, y expr= 1.7*(\thisrow{x} - 1)+1] {data/fig_sine_swarm_lhs.dat};
	
	\addplot[red, only marks, mark=*, mark size = 1pt]
		table [y expr = \thisrow{x}+2, x expr=10*\thisrow{y}] {data/fig_sine_swarm_minimax.dat};

\end{groupplot}
\end{tikzpicture}

\caption{In high dimensional spaces it is critical to incorporate 
structure to construct good experimental designs.
In this example we consider ridge function~\cref{eq:sine} 
and measure the quality of uniform random sampling and Latin hypercube sampling (LHS)
in the Lipschitz matrix metric.
}
\label{fig:sine}
\end{figure}
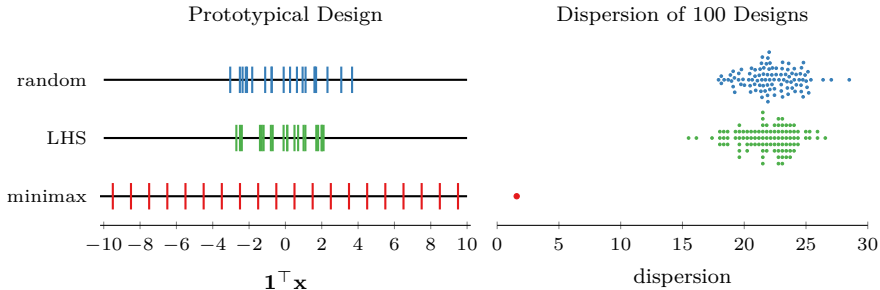

Consider a ridge function
\begin{equation}\label{eq:sine}
	f:[-1,1]^{10} \to \R, \qquad f(\ve x) = \sin(\ve 1^\trans \ve x)
\end{equation}
with a one-dimensional active subspace $\Span \lbrace \ve 1 \rbrace$
and a rank-one Lipschitz matrix $\ma L = \ve 1 \ve 1^\trans$.
If we randomly select points $\lbrace \ve x_j\rbrace_{j=1}^M$
with uniform probability over the domain,
their sum tends to concentrate around the mean (zero)
when projected onto the active subspace
as seen in \cref{fig:sine}.
The same is true with a Latin hypercube design
where points are randomly selected
so that projection onto the coordinate axes
results in evenly spaced points.
However, this is not true when we construct 
a minimax design under the Lipschitz matrix metric.
This design has much lower dispersion
and consequently we can construct much more accurate approximations
according to \cref{thm:error_bound}.

\subsection{Constructing Minimax Designs}
Finding even an approximately optimal minimax design is challenging:
it requires solving a deeply nested optimization problem~\cref{eq:minimax2}
in a high dimensional space
with many local minimizers with large dispersion.
Here we briefly describe a combination of three algorithms
that we use to construct the minimax designs appearing 
in the remainder of this paper.
Each algorithm primarily uses bounded Voronoi vertices under the Lipschitz matrix metric.
All such vertices can be computed using Qhull~\cite{BDH96}
or a random subset by sequential projection~\cite{LC05};
using sequential projection allows us to employ this combination algorithms
spaces of moderate dimension, i.e., $m> 4$. 

To avoid finding a local minimizer with a large dispersion,
the first two algorithm construct a good initial design
that is then refined to local optimality.
First, we construct a maximin coffeehouse design~\cite{Mul01}
where we greedily add new points to the design maximizing their distance
from existing points in the design
\begin{equation}
	\ve x_k \leftarrow \argmax_{\ve x \in \set D} \min_{j=1,\ldots,k-1} \|\ma L(\ve x - \ve x_j)\|_2.
\end{equation}
This optimizer $\ve x_k$ is a bounded Voronoi vertex
that we can (approximately) identify by a finite minimization over (a subset of) these vertices.
Second, we perform a few iterations of block
coordinate descent~\cite{SHSV03} to maximize the minimum pairwise distance
\begin{equation}
	\ve x_k^{\ell+1} \leftarrow \argmax_{\ve x\in \set D}\min_{\substack{j=1,\ldots, M\\ j\ne k}}
		\|\ma L(\ve x - \ve x_j^{\ell})\|_2
\end{equation}
Again, we can (approximately) solve this optimization problem
by a finite minimization over (a subset of) the bounded Voronoi vertices.
Finally this design initializes a variant of Lloyd's algorithm 
to find a locally optimal minimax design~\cite{SD96}.
At each iteration we identify the bounded Voronoi vertices $\ve v_i^\ell$
associated with the design $\lbrace \ve x_j^\ell \rbrace_{j=1}^M$
and then move $\ve x_j^\ell$ to be the circumcenters of its Voronoi cell
\begin{equation}
	\ve x_k^{\ell+1} \leftarrow \argmin_{\ve x\in \set D} 
		\max_{\ve v_i^\ell \in \text{Voronoi cell of }\ve x_k^\ell}
		\| \ma L(\ve x - \ve v_i^\ell)\|_2.
\end{equation} 
If this iteration converges, the design satisfies the local optimality conditions
for a minimax design~\cite[Thm.~4.7]{CB05}.
Note that if the Lipschitz matrix is low-rank
(i.e., rank one, two, or three)
it feasible to compute all the bounded Voronoi vertices
even if the ambient space is high dimensional.

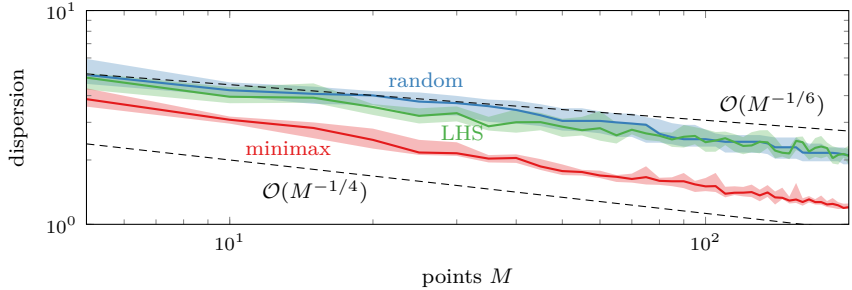
\begin{figure}
\centering
\pgfplotstableread{data/fig_design_rate_random.dat}\random
\pgfplotstableread{data/fig_design_rate_lhs.dat}\lhs
\pgfplotstableread{data/fig_design_rate_minimax.dat}\minimax
\begin{tikzpicture}
\begin{axis}[
		width =0.9\linewidth,
		height = 0.34\linewidth,
		xmin = 5, xmax = 200,
		ymin = 1, ymax = 10,
		ytick = {1,10},
		xlabel = points $M$,
		ylabel = dispersion,
		xmode = log, ymode = log,
	]

	\addplot[blue, thick] table [ x= N, y = p50] {\random}
		node [pos = 0.43, anchor = south, yshift=2pt]{random};

	\addplot[draw=none, name path=random_lb] table [ x=N, y = p25] {\random};
	\addplot[draw=none, name path=random_ub] table [ x=N, y = p75] {\random};
	\addplot[blue, opacity = 0.3] fill between [of=random_lb and random_ub];
	
	\addplot[green, thick] table [ x= N, y = p50] {\lhs}
		node [pos = 0.43, anchor = north, yshift=0pt]{LHS};
	\addplot[draw=none, name path=lhs_lb] table [ x=N, y = p25] {\lhs};
	\addplot[draw=none, name path=lhs_ub] table [ x=N, y = p75] {\lhs};
	\addplot[green, opacity = 0.3] fill between [of=lhs_lb and lhs_ub];
	
	\addplot[red, thick] table [ x= N, y = p50] {\minimax}
		node [pos = 0.25, anchor = north, yshift=-2pt]{minimax};
	\addplot[draw=none, name path=minimax_lb] table [ x=N, y = p25] {\minimax};
	\addplot[draw=none, name path=minimax_ub] table [ x=N, y = p75] {\minimax};
	\addplot[red, opacity = 0.3] fill between [of=minimax_lb and minimax_ub];
	
	\addplot[black, densely dashed, domain = 5:2e2] {4.5*(x/10)^(-1/6)} 
		node [pos = 0.9, anchor = south, rotate=0] {$\order(M^{-1/6})$};
	
	\addplot[black, densely dashed, domain = 5:2e2] {2*(x/10)^(-1/4)} 
		node [pos = 0.3, anchor = north, rotate=0] {$\order(M^{-1/4})$};

\end{axis}
\end{tikzpicture}

\caption{The dispersion of several types of designs
for the OTL circuit function using twenty trials.
The solid line shows the median and the shaded regions the 25th to 75th percentile.
}
\label{fig:design_rate}
\end{figure}
 
It is tempting to think that a minimax design is only necessary
when function evaluations are costly.
However, as \cref{fig:design_rate} illustrates
increasing the number of random and Latin hypercube samples 
decreases dispersion only very slowly at $\order(M^{-1/6})$
whereas we see the minimax design temporally converges
faster and has a substantially reduced dispersion (see discussion in
subsections~\ref{sec:ibc:asymptotic} and~\ref{sec:ibc:nonasymptotic}). 
Hence if we seek a low-dispersion design 
for accurate approximation, integration, or optimization
it is critical to construct an approximate minimax design.

 \section{Parameter Reduction\label{sec:dimension}}
Using the Lipschitz matrix,
we can perform dimension reduction by constructing a ridge approximation
\begin{equation}
	f(\ve x) \approx \widetilde f(\ve x) \coloneqq g(\ma U^\trans \ve x),
	\quad g:\R^n\to \R, \quad \ma U \in \R^{m\times n}.
\end{equation}
Since $g$ is a ridge function, $g$ has at most a rank-$n$ Lipschitz matrix
(\cref{thm:ridge})
and complexity of tasks now scales with the ridge dimension $n$
and no longer the dimension of the ambient space $m$.
There are many ways to construct ridge approximations;
e.g., 
picking $\ma U$ from the dominant eigenvectors of the MeGO matrix
and then fitting a polynomial~\cite{Con15}
or directly picking $\ma U$ and polynomial $g$ via optimization~\cite{CEHW17, HC18}.
Here we construct a ridge approximation using the Lipschitz matrix
in light of \cref{thm:error_bound}.

\subsection{Building a Ridge Approximation\label{sec:dimension:build}}
Our choice of the ridge subspace is clear from \cref{thm:error_bound}:
the leading $n$ eigenvectors of the Lipschitz matrix form $\ma U\in \R^{m\times n}$,
same as when using the MeGO matrix.
To construct a ridge approximation 
satisfying \cref{thm:error_bound}
we need to choose a ridge function $\widetilde f\in \set L(\set D, \ma L\ma U\ma U^\trans)$.
Due to the constraints its Lipschitz matrix, $\widetilde f$ may not interpolate point queries;
i.e., $\widetilde f(\ve x_j) \ne f(\ve x_j)$.
Hence to a compatible ridge approximation $\widetilde f$
we first solve an optimization problem for 
function values $\widetilde y_j = \widetilde f(\ve x_j)$
minimizing the $\ell_1$ error (or equivalently any other convex norm) 
\begin{equation}\label{eq:Lipschitz_profile}
	\tve y = \argmin_{\ve y\in \R^M}  \sum_{j=1}^M | y_j - f(\ve x_j)|, \quad
	\text{such that} \quad |y_j - y_k| \le \|\ma L\ma U\ma U^\trans(\ve x_j - \ve x_k)\|_2.
\end{equation}
Then these function values 
define point queries $\tset P_M = \lbrace \ve x_j,  \widetilde y_j\rbrace_{j=1}^M$
and define $\widetilde f$ to be the central approximation of this data,
cf.~\cref{eq:central_approx}:
\begin{equation}\label{eq:Lipschitz_approx}
	\widetilde f(\ve x) =
		\frac12 \left[
			\fmin(\ve x; \ma L \ma U\ma U^\trans, \tset P_M)
			+
			\fmax(\ve x; \ma L \ma U\ma U^\trans, \tset P_M)
		\right].
\end{equation}
As $\widetilde f$ approximates $f$ with a maximum error $\epsilon = \max_{j} |y_j - f(\ve x_j)|$,
the uncertainty associated with $\widetilde f$ is, cf.~\cref{eq:uncertainty_interval},
\begin{equation}\label{eq:uncertainty_set_epsilon}
	\set U_\epsilon(\ve x; \ma L\ma U\ma U^\trans, \tset P_M) 
		= [\fmin(\ve x; \ma L\ma U\ma U^\trans, \tset P_M) - \epsilon \ ,
		 \ \fmax(\ve x; \ma L\ma U \ma U^\trans, \tset P_M) + \epsilon].
\end{equation}

\subsection{Typical Workflow\label{sec:dimension:workflow}}
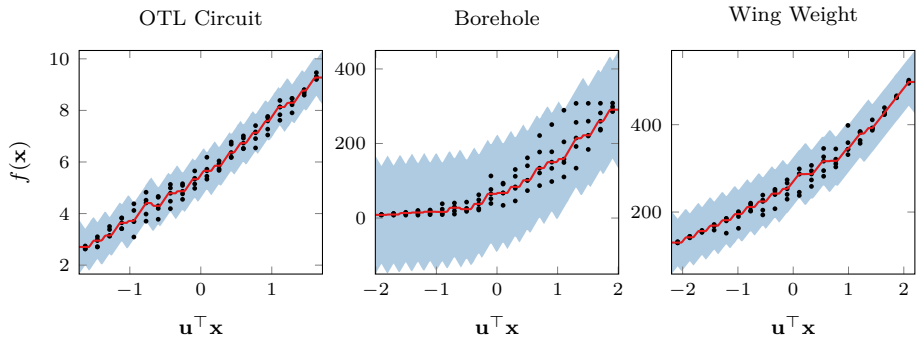
\begin{figure}
\centering
\begin{tikzpicture}
\begin{groupplot}[
		group style = {group size = 3 by 1, horizontal sep = 20pt},
		height = 0.35\linewidth,
		width = 0.37\linewidth,
		xlabel = $\ve u^\trans \ve x$,
		clip mode=individual,
	]
	\nextgroupplot[
		title = OTL Circuit,
		ylabel = $f(\ve x)$,
	]
	\pgfplotstableread{data/fig_ridge2_circuit_curves.dat}\curves
	\pgfplotstableread{data/fig_ridge2_circuit_points.dat}\points
	\pgfplotstableread{data/fig_ridge2_circuit_envelope.dat}\envelope

	\addplot[draw=none, name path = circuit_lb] table [x=x, y=lb] {\curves};
	\addplot[draw=none, name path = circuit_ub] table [x=x, y=ub] {\curves};
	\addplot[blue, opacity = 0.4] fill between [of=circuit_lb and circuit_ub];

	\addplot[thick, black, only marks, mark = *, mark size =0.5pt] table [x=x, y = y]{\points};
	
	\addplot[thick, red] table [x=x, y = c]{\curves};
	
	\nextgroupplot[
		title = Borehole,
	]
	\pgfplotstableread{data/fig_ridge2_borehole_curves.dat}\curves
	\pgfplotstableread{data/fig_ridge2_borehole_points.dat}\points

	\addplot[draw=none, name path = borehole_lb] table [x=x, y=lb] {\curves};
	\addplot[draw=none, name path = borehole_ub] table [x=x, y=ub] {\curves};
	\addplot[blue, opacity = 0.4] fill between [of=borehole_lb and borehole_ub];

	\addplot[thick, black, only marks, mark = *, mark size =0.5pt] table [x=x, y = y]{\points};
	
	\addplot[thick, red] table [x=x, y = c]{\curves};
	
	\nextgroupplot[
		title = Wing Weight,
	]
	\pgfplotstableread{data/fig_ridge2_wing_curves.dat}\curves
	\pgfplotstableread{data/fig_ridge2_wing_points.dat}\points

	\addplot[draw=none, name path = wing_lb] table [x=x, y=lb] {\curves};
	\addplot[draw=none, name path = wing_ub] table [x=x, y=ub] {\curves};
	\addplot[blue, opacity = 0.4] fill between [of=wing_lb and wing_ub];

	\addplot[thick, black, only marks, mark = *, mark size =0.5pt] table [x=x, y = y]{\points};
	
	\addplot[thick, red] table [x=x, y = c]{\curves};

\end{groupplot}
\end{tikzpicture}
\caption{Three ridge approximations as described in \cref{sec:dimension:workflow}.
In each plot the $x$-coordinate shows the projection of the domain 
onto the one-dimensional active subspace.
The red curve denotes the central approximation,
the black dots the data used to construct the approximation,
and the shaded region the uncertainty.}
\label{fig:ridge}
\end{figure}

A typical workflow for using the Lipschitz matrix for dimension reduction
starts with querying the gradient on the corners of the domain
and using this data to estimate the Lipschitz matrix.
Then, picking the active subspace to be the span of the leading eigenvector $\ve U\in \R^{m\times 1}$ of $\ma L$,
we construct a 20-point minimax design under this metric $\ma L\ve U \ve U^\trans$,
$\lbrace \ve x_j \rbrace_{j=1}^{20}$.
To better estimate the variability at each point,
we take 5 maximin coffeehouse samples 
on each domain $\set D_j = \lbrace \ve x\in \set D: \ma U^\trans \ve x = \ma U^\trans \ve x_j\rbrace$
in the Lipschitz metric $\ma L$.
Finally, we construct the ridge approximation as described in the previous section.
\Cref{fig:ridge} shows three examples of this workflow
and the associated uncertainty.

 \section{Uncertainty Quantification\label{sec:uq}}
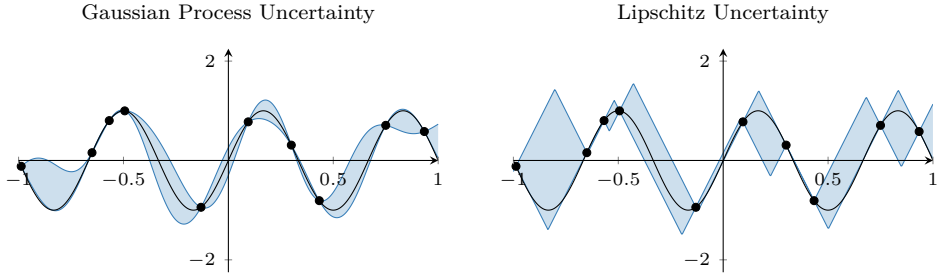
\begin{figure}[t]
\centering
\begin{tikzpicture}
\begin{groupplot}[
	group style = {group size = 2 by 1},
	width = 0.55\textwidth,	
	height = 0.35\textwidth,
	axis lines = middle,
	ymin = -2.25,ymax = 2.25,
	]
	\nextgroupplot[title=Gaussian Process Uncertainty]
	
	\addplot[blue, name path=gpr_lb] table [x=x, y=gpr_lb] {data/fig_gp.dat};	
	\addplot[blue, name path=gpr_ub] table [x=x, y=gpr_ub] {data/fig_gp.dat};	
	\addplot[blue, opacity = 0.25] fill between [of=gpr_lb and gpr_ub];

	\addplot[black, mark=*, only marks, mark size=1.5] table [x=x, y=fx] {data/fig_gp_data.dat};
	\addplot[black] table [x=x, y=fx] {data/fig_gp.dat};

	\nextgroupplot[title=Lipschitz Uncertainty]
	
	\addplot[blue, name path=lip_lb] table [x=x, y=lip_lb] {data/fig_gp.dat};	
	\addplot[blue, name path=lip_ub] table [x=x, y=lip_ub] {data/fig_gp.dat};	
	\addplot[blue, opacity = 0.25] fill between [of=lip_lb and lip_ub];
	
	\addplot[black, mark=*, only marks, mark size=1.5] table [x=x, y=fx] {data/fig_gp_data.dat};
	\addplot[black] table [x=x, y=fx] {data/fig_gp.dat};

\end{groupplot}
\end{tikzpicture}

\caption[A comparison of Gaussian process and Lipschitz uncertainty]{	A comparison of Gaussian Process and the Lipschitz notions of uncertainty.
	Here we use ten samples of $f(x)=\sin(3\pi x)$ on $\set D= [-1,1]$ denoted by dots,
	with $f$ shown by the black line.
	In both plots, the shaded area denotes the uncertainty estimate.
	The left plot shows the Gaussian process with
	zero-mean $m(\ve x)=0$ and
	squared exponential covariance $k(\ve x, \ve x') = \exp(-\frac12 \| \ell(\ve x - \ve x')\|_2^2)$
	where $\ell=0.157$ has been chosen to maximize the marginal likelihood~\cite[subsec.~2.7.1]{RW06};
	the probability threshold $\delta$ corresponds to one-standard deviation.
	The right plot shows the Lipschitz uncertainty with Lipschitz constant $L=8.39$ estimated from the samples.
}
\label{fig:gp}
\end{figure}

When working with expensive deterministic computer simulations
it is often necessary to employ an approximation of certain quantities of interest, 
called a \emph{response surface} or a \emph{surrogate}.
Supposing we have constructed this approximation using samples of $f$,
it is natural ask: what are the range of possible values our approximation could take away from these samples?
This is often called \emph{uncertainty} in this setting.
Gaussian processes provide one approach to define an uncertainty~\cite[sec.~2.2]{RW06}
and the Lipschitz constant provides another~\cite{RS15}.
These two techniques are based on different assumptions
and yield different results as illustrated in \cref{fig:gp}.
The Gaussian process approach assumes that 
$f$ is a Gaussian process conditioned on the measurements $y_j = f(\ve x_j)$,
maximizes the likelihood of a parameterized covariance kernel based on observations,
and then computes the probability of observing $f(\ve x)$;
the uncertainty is visualized as all function values above some probability threshold.
In contrast, the Lipschitz approach assumes measurements $y_j= f(\ve x_j)$ comes from a Lipschitz function,
chooses the smallest Lipschitz constant that consistent with this data,
and defines uncertainty as the range of possible function values
consistent with this Lipschitz constant and data;
namely, the \emph{uncertainty set} $\set U$~\cref{eq:uncertainty_set}.
Replacing the Lipschitz constant
by the Lipschitz matrix
allows us to reduce uncertainty.

\subsection{Decreasing Uncertainty}
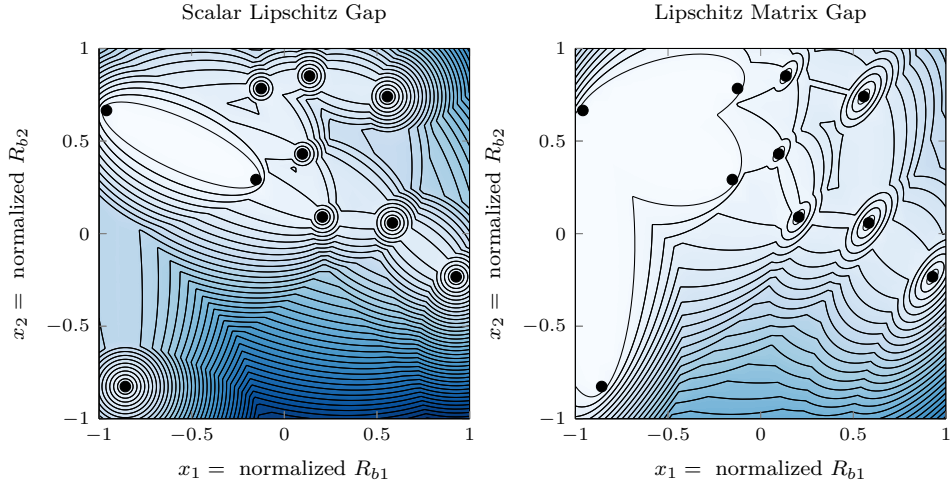
\begin{figure}
\centering 
\begin{tikzpicture}
\begin{groupplot}[
	group style = {group size = 2 by 1, horizontal sep=4em},
	width = 0.5\textwidth,	
	height = 0.5\textwidth,
	axis equal,
	view={0}{90},
	colormap name=Blu9cont,
	point meta min = 0,
	point meta max = 3.71,
	xlabel = {$x_1= \text{ normalized } R_{b1}$},
	ylabel = {$x_2= \text{ normalized } R_{b2}$},
	]
	\nextgroupplot[title = Scalar Lipschitz Gap]
	\addplot3[contour filled={number = 100, labels = {false}},
		mesh/rows=20, mesh/cols=20, mesh/check=false,]
		table {data/fig_gap_scalar_otl_circuit.dat};
	\addplot[black, only marks, mark=*]
		table {data/fig_gap_scalar_otl_circuit_samples.dat};
	
	\addplot[contour prepared ={ draw color=black,labels=false}, contour prepared format=matlab] 
		table {data/fig_gap_scalar_otl_circuit_contour.dat};
	
	\nextgroupplot[title = Lipschitz Matrix Gap]
	\addplot3[contour filled={number = 100, labels = {false}},
		mesh/rows=20, mesh/cols=20, mesh/check=false,]
		table {data/fig_gap_matrix_otl_circuit.dat};
	\addplot[black, only marks, mark=*]
		table {data/fig_gap_matrix_otl_circuit_samples.dat};
	
	\addplot[contour prepared ={ draw color=black,labels=false}, contour prepared format=matlab] 
		table {data/fig_gap_matrix_otl_circuit_contour.dat};

\end{groupplot}
\end{tikzpicture}
\caption[Lipschitz Bound Gap]{A plot of the gap 
between the upper and lower Lipschitz bounds.
Here we use function values taken the OTL circuit test problem using the first two parameters,
holding the remainder at their central values;
additionally the domain has been normalized to $[-1,1]^2$.
In each plot, the contours show an increase of $0.1$ in the gap.
Here the uncertainty set, the Lipschitz constant, and Lipschitz matrix
were all estimated using the ten samples marked with dots.
Here $L=2.38$
and 
$\ma L = \begin{bsmallmatrix} 1.6 & 0 \\ -1.3 & 1.6\end{bsmallmatrix}$.
}
\label{fig:gap}
\end{figure}

As the Lipschitz matrix more accurately encodes
variability in the function, 
it can reduces uncertainty 
compared to the scalar Lipschitz constant.
\Cref{fig:gap} illustrates the improvement of the 
Lipschitz matrix over the Lipschitz constant on a 
two dimensional example.
\Cref{tab:sample} demonstrates the same reduction in uncertainty
occurs in higher dimensional examples.
This example shows that uncertainty is reduced the most
when Lipschitz matrix is used both when
constructing the minimax design
and evaluating the uncertainty set.

\begin{table}
\caption{Estimated maximum pointwise uncertainty over the domain
using a 100-point minimax design with either isotropic or Lipschitz matrix metric.
The maximum uncertainty is computed by extensively sampling the domain
and measuring the diameter the Lipschitz constant or Lipschitz matrix uncertainty set. 
To aid comparison, each function's range has been normalized to $[0,1]$. 
}
\label{tab:sample}
\centering

\pgfplotstablevertcat{\Lmat}{data/tab_sample_golinski_matrix_lipschitz_uncertainty.dat}
\pgfplotstablevertcat{\Lmat}{data/tab_sample_otl_matrix_lipschitz_uncertainty.dat}
\pgfplotstablevertcat{\Lmat}{data/tab_sample_piston_matrix_lipschitz_uncertainty.dat}
\pgfplotstablevertcat{\Lmat}{data/tab_sample_borehole_matrix_lipschitz_uncertainty.dat}
\pgfplotstablevertcat{\Lmat}{data/tab_sample_wing_matrix_lipschitz_uncertainty.dat}

\pgfplotstablecreatecol[expr = \thisrow{p100}] {Lmat}{\Lmat}

\pgfplotstablevertcat{\Liso}{data/tab_sample_golinski_matrix_isotropic_uncertainty.dat}
\pgfplotstablevertcat{\Liso}{data/tab_sample_otl_matrix_isotropic_uncertainty.dat}
\pgfplotstablevertcat{\Liso}{data/tab_sample_piston_matrix_isotropic_uncertainty.dat}
\pgfplotstablevertcat{\Liso}{data/tab_sample_borehole_matrix_isotropic_uncertainty.dat}
\pgfplotstablevertcat{\Liso}{data/tab_sample_wing_matrix_isotropic_uncertainty.dat}

\pgfplotstablecreatecol[expr = \thisrow{p100}] {Liso}{\Liso}
\pgfplotstablecreatecol[copy column from table={\Liso}{Liso}] {Liso} {\Lmat}

\pgfplotstablevertcat{\Lcon}{data/tab_sample_golinski_scalar_isotropic_uncertainty.dat}
\pgfplotstablevertcat{\Lcon}{data/tab_sample_otl_scalar_isotropic_uncertainty.dat}
\pgfplotstablevertcat{\Lcon}{data/tab_sample_piston_scalar_isotropic_uncertainty.dat}
\pgfplotstablevertcat{\Lcon}{data/tab_sample_borehole_scalar_isotropic_uncertainty.dat}
\pgfplotstablevertcat{\Lcon}{data/tab_sample_wing_scalar_isotropic_uncertainty.dat}

\pgfplotstablecreatecol[expr = \thisrow{p100}] {Lcon}{\Lcon}
\pgfplotstablecreatecol[copy column from table={\Lcon}{Lcon}] {Lcon} {\Lmat}

\pgfplotstablevertcat{\Lrange}{data/tab_sample_golinski_scalar_isotropic_uncertainty.dat}
\pgfplotstablevertcat{\Lrange}{data/tab_sample_otl_scalar_isotropic_uncertainty.dat}
\pgfplotstablevertcat{\Lrange}{data/tab_sample_piston_scalar_isotropic_uncertainty.dat}
\pgfplotstablevertcat{\Lrange}{data/tab_sample_borehole_scalar_isotropic_uncertainty.dat}
\pgfplotstablevertcat{\Lrange}{data/tab_sample_wing_scalar_isotropic_uncertainty.dat}

\pgfplotstablecreatecol[expr = \thisrow{range}] {range}{\Lrange}

\pgfplotstablecreatecol[copy column from table={\Lrange}{range}] {range} {\Lmat}

\pgfplotstabletypeset[
	every head row/.style={before row=\toprule,after row=\midrule},
    every last row/.style={after row=\bottomrule},
	create on use/name/.style={
        create col/set list={
			\hspace*{-2pt}Golinski volume~\cite{Gol70},
			OTL circult~\cite{BS07a}, 
			piston~\cite{KZ98},
			borehole~\cite{HG83}, 
			wing weight~\cite{FSK08},
		},
    },
	create on use/dimension/.style={
        create col/set list={
			6, 6, 7, 8, 10,
		},
	},
	create on use/ratio/.style={
      create col/expr={\thisrow{Lcon}/\thisrow{Lmat}}},
	columns/name/.style = {column name = test problem, string type, column type = {r}, },
	columns/dimension/.style = {int detect, column name = {\hspace*{-2pt}dim.\hspace*{-2pt}}  },
	create on use/Lconnorm/.style={
      create col/expr={\thisrow{Lcon}/\thisrow{range}}},
	create on use/Lmatnorm/.style={
      create col/expr={\thisrow{Lmat}/\thisrow{range}}},
	create on use/Lisonorm/.style={
      create col/expr={\thisrow{Liso}/\thisrow{range}}},
	columns/Lconnorm/.style = {column name = \makecell[b]{isotropic \\ Lip. const. \\ uncertainty}, zerofill, precision=2, dec sep align, fixed},
	columns/Lmatnorm/.style = {column name = \makecell[b]{Lipschitz \\ Lip. matrix \\ uncertainty}, zerofill, precision=2, dec sep align, fixed},
	columns/Lisonorm/.style = {column name = \makecell[b]{isotropic \\ Lip. matrix \\ uncertainty}, zerofill, precision=2, dec sep align, fixed},
	columns/ratio/.style = {column name = ratio, zerofill, precision=2, dec sep align, fixed},
	columns = {name, dimension, Lconnorm, Lisonorm, Lmatnorm, ratio},
]{\Lmat}

\end{table}

\subsection{Visualizing Uncertainty on Shadow Plots}
\begin{figure}
\centering
\pgfplotstableread{data/fig_shadow_iso_scalar_shadow_uncertainty.dat}\uncertainA
\pgfplotstableread{data/fig_shadow_iso_matrix_shadow_uncertainty.dat}\uncertainB
\pgfplotstableread{data/fig_shadow_lip_matrix_shadow_uncertainty.dat}\uncertainC
\begin{tikzpicture}
\begin{groupplot}[
		group style = {group size = 3 by 1, horizontal sep = 20pt},
		height = 0.35\linewidth,
		width = 0.38\linewidth,
		xlabel = $\ve u^\trans \ve x$,
		clip mode=individual,
		xmin = -2,
		xmax = 2,
		ymin = -2,
		ymax = 12,
		title style={align=center}, 
	]

	\nextgroupplot[title = {Isotropic Design \\ Lipschitz Constant},
		ylabel = $f(\ve x)$, 
		]	
		
	\addplot[only marks, mark=*, mark size =0.5pt]
		table [x=y, y=fX] {data/fig_shadow_iso_scalar.dat};
	
	\addplot[green, dotted, name path = lb] table [x=y, y=lb] {\uncertainA};
	\addplot[green, dotted, name path = ub] table [x=y, y=ub] {\uncertainA};
	\addplot[green, opacity = 0.4] fill between [of=lb and ub];
	
	\addplot[green, dashed] table [x=y, y=lb] {\uncertainB};
	\addplot[green, dashed] table [x=y, y=ub] {\uncertainB};
	
	\addplot[green] table [x=y, y=lb] {\uncertainC};
	\addplot[green] table [x=y, y=ub] {\uncertainC};
	
	\nextgroupplot[title = {Isotropic Design \\ Lipschitz Matrix} 
		]	
		
	\addplot[only marks, mark=*, mark size =0.5pt]
		table [x=y, y=fX] {data/fig_shadow_iso_matrix.dat};
	
	\addplot[green, dashed, name path = lb] table [x=y, y=lb] {\uncertainB};
	\addplot[green, dashed, name path = ub] table [x=y, y=ub] {\uncertainB};
	\addplot[green, opacity = 0.4] fill between [of=lb and ub];
	
	\addplot[green, dotted] table [x=y, y=lb] {\uncertainA};
	\addplot[green, dotted] table [x=y, y=ub] {\uncertainA};
	
	\addplot[green] table [x=y, y=lb] {\uncertainC};
	\addplot[green] table [x=y, y=ub] {\uncertainC};
	
	\nextgroupplot[title = {Lipschitz Matrix Design \\ Lipschitz Matrix} 
		]	
		
	\addplot[only marks, mark=*, mark size =0.5pt]
		table [x=y, y=fX] {data/fig_shadow_lip_matrix.dat};
	
	\addplot[green, name path = lb] table [x=y, y=lb] {\uncertainC};
	\addplot[green, name path = ub] table [x=y, y=ub] {\uncertainC};
	\addplot[green, opacity = 0.4] fill between [of=lb and ub];

	\addplot[green, dotted] table [x=y, y=lb] {\uncertainA};
	\addplot[green, dotted] table [x=y, y=ub] {\uncertainA};
	
	\addplot[green, dashed] table [x=y, y=lb] {\uncertainB};
	\addplot[green, dashed] table [x=y, y=ub] {\uncertainB};

\end{groupplot}
\end{tikzpicture}

\caption{Projected shadow plots for the OTL Circuit problem
along the leading eigenvector of the Lipschitz matrix.
Black dots denote the design projected onto this active subspace $\Span \ve u$,
the dotted, dashed, and solid lines denote
the combination of an isotopic minimax design and the scalar Lipschitz constant uncertainty, 
an isotropic minimax design and the Lipschitz matrix uncertainty,
and a minimax design in the Lipschitz matrix metric and the Lipschitz matrix uncertainty.
The shaded region on each plot represents the projected uncertainty.}
\label{fig:shadow}
\end{figure}
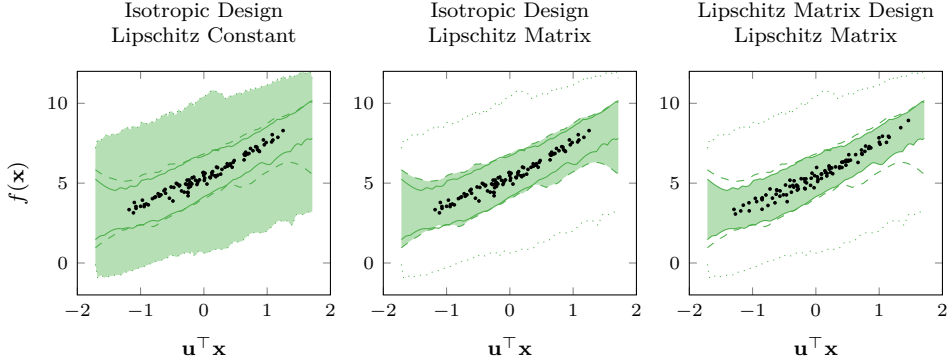

Shadow plots are an important tool for visualizing functions with a high-dimensional domain.
To include the uncertainty in shadow plots,
we generalize the uncertainty set~\cref{eq:uncertainty_set}
to take set-valued inputs
\begin{equation}\label{eq:uncertainty_set_set}
	\begin{split}
		\set U(\set S; \ma L, \set P_M) 
			&\coloneqq \bigcup_{\ve x \in \set S} \set U(\ve x; \ma L, \set P_M)\\
			& = \left[ \min_{\ve x\in \set S} \max_j y_j - \|\ma L(\ve x - \ve x_j)\|_2 \ , \ 
				\max_{\ve x \in \set S} \min_j y_j + \|\ma L(\ve x - \ve x_j)\|_2 \right].
	\end{split}
\end{equation}
In our implementation, 
we identify these lower and upper bounds by solving a sequential linear program~\cite{OW69}.
Then to include the uncertainty in a one-dimensional shadow plot
we evaluate the uncertainty set for 
$\set S_\alpha = \lbrace \ve x\in \set D: \ve u^\trans \ve x = \alpha\rbrace$
for multiple $\alpha$ in the projection of the domain onto $\ve u$.
\Cref{fig:shadow} provides an example of this projected uncertainty 
using both the Lipschitz constant and Lipschitz matrix
as well as with an isotropic minimax design or a Lipschitz matrix minimax design.
Note the uncertainty in \cref{fig:shadow} is the projection of a high-dimensional
uncertainty~\cref{eq:uncertainty_set_set}
rather than an intrinsically low-dimensional uncertainty~\cref{eq:uncertainty_set_epsilon}
in \cref{fig:ridge}.

  \section{$\epsilon$-Lipschitz\label{sec:epsilon}}
In this section we explore three applications of the $\epsilon$-Lipschitz matrix
exploiting its ability to ignore a small changes in the function.

\subsection{Computational Noise\label{sec:epsilon:noise}}
Many functions appearing in computational science and engineering
often have \emph{computational noise}~\cite{MW11}---a phenomenon emerging from many factors
including convergence tolerances and mesh discretizations.
With the addition of computational noise, 
functions that are ideally Lipschitz continuous
can become discontinuous.
The $\epsilon$-Lipschitz matrix enables us to perform 
dimension reduction for these discontinuous functions
and remove the influence of noise.
As an example, consider the partial trace function of Mor\'e and Wild~\cite[Sec.~1]{MW11}:
the sum of the first five smallest eigenvalues of the parameterized matrix 
\begin{equation}
	f: [-1,1]^2 \to \R, \qquad f(\ve x) = \sum_{i=0}^4 \lambda_{700-i}(
		\ma A + 2.5(x_1+1) \ve e_1\ve e_1^\trans + 2.5(x_2+1) \ve e_2 \ve e_2^\trans
		)
\end{equation}
where $\ma A\in \R^{700\times 700}$ is the \verb|Trefethen_700| sparse matrix~\cite{DH11},
and $\lambda_i$ denotes the $i$th eigenvalue in decreasing order.
If we accurately evaluate this function
it is approximately linear on the domain
and yields an approximately low-rank Lipschitz matrix.
To introduce computational noise,
we use a very loose relative accuracy termination criteria of $0.1$
when computing these eigenvalues using ARPACK~\cite{LSY98}.
As evidenced in \cref{tab:noise},
the introduction of noise pollutes the estimate of the Lipschitz matrix,
yielding a large, full-rank Lipschitz matrix.
If we estimate an $\epsilon$-Lipschitz matrix instead with $\epsilon=6.47$ 
(the largest mismatch between accurate and noisy evaluations of this function)
we are able to recover an accurate, approximately low-rank Lipschitz matrix.
In this example we use samples on a $10\times 10$ grid
and initialize the Lanczos iteration using the ones vector.

\begin{table}
\caption{Computational noise can lead to inaccurate estimates of the Lipschitz matrix.
Using the $\epsilon$-Lipschitz matrix instead, 
we can recover an accurate Lipschitz matrix.
These Lipschitz matrices correspond to the partial trace example
described in \cref{sec:epsilon:noise}.
}
\label{tab:noise}

\centering
\begin{tabular}{p{0.28\linewidth} p{0.28\linewidth} p{0.32\linewidth}}
\toprule
Exact data & Noisy data & Noisy data with $\epsilon$-Lipschitz \\
\midrule 
$\ma L_{\phantom{1}} = \begin{bmatrix}
1.89 &  1.63 \\
1.63 &  1.84
\end{bmatrix}$ 
&
$ \ma L_{\phantom{1}} = \begin{bmatrix}
	31.6 &  0.118 \\
	0.118 & 31.5 \\
\end{bmatrix}$
&
$\ma L_\epsilon=\begin{bmatrix}
1.84 & 1.66 \\
1.66 & 1.79 
\end{bmatrix}$
\\[10pt]
$\lambda_1 = 3.50, \ \lambda_2 = 0.233$ &
$\lambda_1 = 31.7, \ \lambda_2 = 31.4$ & 
$\lambda_1 = 3.48, \ \lambda_2 = 0.151$ \\
\bottomrule
\end{tabular}
\end{table}

\subsection{Distracting Oscillations\label{sec:epsilon:oscillations}}
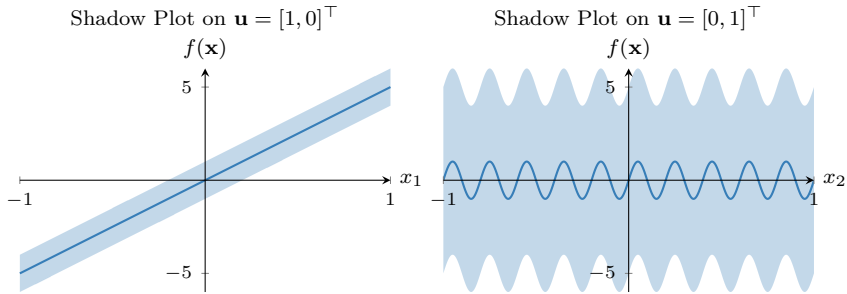
\begin{figure}
\centering
\begin{tikzpicture}
\begin{groupplot}[
	group style = {group size = 2 by 1, 
		horizontal sep = 20pt,
		vertical sep = 10pt},
	width = 0.5\textwidth,	
	height = 0.35\textwidth,
	ymin = -6,
	ymax = 6,
	xmin = -1,
	xmax = 1,
	axis lines = middle,
	xtick = {-1,1},
	ylabel = {$f(\ve x)$},
	x label style={at={(axis cs:1,0)},anchor=west},
	y label style={at={(axis cs:0,6)},anchor=south},
	title style={yshift=5pt},
	]

	\nextgroupplot[title = {Shadow Plot on $\ve u = [1,0]^\trans$},
		xlabel = {$x_1$},
	]
	
	\addplot[colorbrewerA2, thick, domain = -1:1] {5*x};
	\addplot[draw=none, domain = -1:1, name path=e1lb] {5*x-1};
	\addplot[draw=none, domain = -1:1, name path=e1ub] {5*x+1};
	\addplot[colorbrewerA2, opacity = 0.3] fill between [of=e1lb and e1ub];
	
	\nextgroupplot[title = {Shadow Plot on $\ve u = [0,1]^\trans$},
		samples = 1000,
		xlabel = {$x_2$},
	]
	
	\addplot[colorbrewerA2, thick, domain = -1:1] {sin(10*180*x)};
	\addplot[draw=none, domain = -1:1, name path=e2lb] {sin(10*180*x)-5};
	\addplot[draw=none, domain = -1:1, name path=e2ub] {sin(10*180*x)+5};
	\addplot[colorbrewerA2, opacity = 0.3] fill between [of=e2lb and e2ub];

\end{groupplot}
\end{tikzpicture}

\caption{Two shadow plots of the corrugated roof function~\cref{eq:roof}.
The solid line denotes the mean value,
the shaded area denotes the area of possible oscillation.
}
\label{fig:roof}
\end{figure}

The active subspace identified using a Lipschitz matrix 
can yield counterintuitive results when a function
has high frequency, low amplitude oscillations.
Consider the ``corrugated roof'' function~\cite[eq.~(26)]{CEHW17}:
\begin{equation}\label{eq:roof}
	f: [-1,1]^2 \to \R, \qquad f(\ve x) = 5 x_1 + \sin(10 \pi x_2). 
\end{equation}
Most of the variation in this function is due to $x_1$,
but there is a high frequency oscillation in $x_2$.
As $f$ is additive in the two coordinates,
we can identify a Lipschitz matrix 
$\ma L = \begin{bsmallmatrix}
	5 & 0 \\ 0 &10 \pi 
\end{bsmallmatrix}$
by considering the largest derivative in each coordinate.
Then applying~\cref{thm:error_bound},
we would choose the active subspace $\Span \lbrace \ve e_2 \rbrace$.
However, as illustrated in \cref{fig:roof},
this yields much higher uncertainty than using $\Span \lbrace \ve e_1 \rbrace$
because the oscillations in $x_2$ have less influence on the output
than the linear term in $x_1$.
We can remove the influence of these oscillations by using the $\epsilon$-Lipschitz matrix.
Taking $\epsilon=2$, we can ignore the sine term
and identify an $\epsilon$-Lipschitz matrix
$\ma L_2 = \begin{bsmallmatrix}
	5 & 0 \\ 0 & 0 
\end{bsmallmatrix}$
With this choice, we choose the active subspace $\Span \lbrace \ve e_1 \rbrace$.
Although this has increased the uncertainty at each point,
it yields better active subspace which reduces overall uncertainty.

\subsection{Dimension Reduction\label{sec:epsilon:dimension}}
\begin{figure}

\centering

\begin{tikzpicture}
\pgfplotstableread{data/fig_epsilon_rank.dat}\data;

\begin{groupplot}[
	group style = {group size = 1 by 2, vertical sep =8pt},
	xmin = 0, xmax = 7,
	width=0.9\linewidth,
	height = 0.2\linewidth,
	]
	\nextgroupplot[
		ymin = 0, ymax = 4,
		height = 0.25\linewidth,
		xticklabels = {,,},
		ylabel = $\|\ma L_\epsilon\|_\fro$,
	]
	\addplot[red, thick, mark=*, mark size=1pt] table [x=epsilon, y = obj] {\data};
	
	\nextgroupplot[
		ymin = 0, ymax = 6.5,
		height = 0.2\linewidth,
		xlabel = $\epsilon$,
		ylabel = rank $\ma L_\epsilon$,
	]
	\addplot[red, thick, mark=*, mark size=1pt] table [x=epsilon, y = rank] {\data};
	
	\addplot[black, opacity = 0.3] coordinates {(0,1) (20,1)}; 
	\addplot[black, opacity = 0.3] coordinates {(0,2) (20,2)}; 
	\addplot[black, opacity = 0.3] coordinates {(0,3) (20,3)}; 
	\addplot[black, opacity = 0.3] coordinates {(0,4) (20,4)}; 
	\addplot[black, opacity = 0.3] coordinates {(0,5) (20,5)}; 
	\addplot[black, opacity = 0.3] coordinates {(0,6) (20,6)}; 
\end{groupplot}

\end{tikzpicture}
\caption{With increasing $\epsilon$, we can identify low-rank $\epsilon$-Lipschitz matrices.
This example shows the OTL Cirucit function
where the $\epsilon$-Lipschitz matrix is estimated using 200 maximin samples 
chosen using an approximate Lipschitz matrix for this function.}
\label{fig:epsilon_rank}
\end{figure}
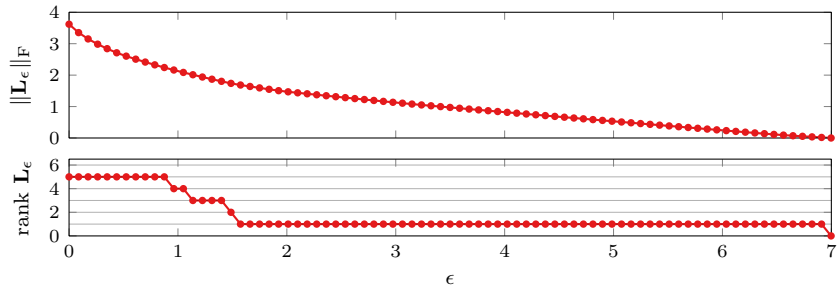

Just as we can use an $\epsilon$-Lipschitz matrix  
to remove undesirable oscillations,
we can also use it for dimension reduction.
If there is a subspace $\set U$ on which a function varies by less than $\epsilon$,
then there is an $\epsilon$-Lipschitz matrix
that has a nullspace of $\dim(\set U)$.
We have no guarantee that by minimizing the Frobenius norm 
with the given function evaluations
we will identify this low-rank $\epsilon$-Lipschitz matrix.
However \Cref{fig:epsilon_rank} illustrates
that we can sometimes find a low-rank matrix.
In this example we note that by ignoring about $20\%$ of the variation,
we can identify a rank-1 $\epsilon$-Lipschitz matrix.

  \section*{Acknowledgements}
The authors would like to thank Akil Narayan and Drew Kouri
for their comments that helped improve this manuscript.
 
\bibliographystyle{siamplain}
\bibliography{abbrevjournals,master,lipschitz,pauls-refs}

\end{document}